\documentclass[sn-mathphys,Numbered]{sn-jnl}

\usepackage[utf8]{inputenc}
\usepackage[T1]{fontenc}
\usepackage{lmodern}
\usepackage{microtype}

\usepackage{amsmath,amssymb,amsfonts,amsthm,mathrsfs}
\usepackage{graphicx}
\usepackage{subcaption}
\usepackage{booktabs}
\usepackage{multirow}
\usepackage{tabularx}
\usepackage{threeparttable}
\usepackage{adjustbox}
\usepackage{enumitem}
\usepackage{pifont}
\usepackage{anyfontsize}
\usepackage[title]{appendix}
\usepackage{xcolor}
\usepackage{placeins}
\usepackage{textcomp}
\usepackage{manyfoot}
\usepackage{listings}
\usepackage[linesnumbered,ruled,vlined]{algorithm2e}
\usepackage{doi}
\usepackage{mdframed}
\usepackage{comment}
\usepackage{hyperref}
\usepackage{cleveref}
\usepackage{url}
\usepackage{nicefrac}

\SetKw{KwIn}{in}
\SetKw{KwOr}{or}
\SetKw{KwAnd}{and}

\SetKwFunction{Function}{\textcolor{awesome}{function}}
\newcommand{\assign}{\leftarrow}

\newcommand{\FuncCall}[2]{\texttt{\bfseries #1(#2)}}
\SetKwProg{Function}{function}{}{}

\DeclareMathOperator*{\argmax}{arg\,max}
\DeclareMathOperator*{\argmin}{arg\,min}
\DeclareMathOperator*{\argsort}{arg\,sort}

\newtheorem{theorem}{Theorem}
\newtheorem{proposition}[theorem]{Proposition}

\newtheorem{definition}{Definition}
\newtheorem{property}{Property}

\begin{document}
\title[Article Title]{An Efficient Global Optimization Algorithm with Adaptive Estimates of the Local Lipschitz Constants}


\author*{\fnm{Danny} \sur{D'Agostino}}
\email{dannydag@nus.edu.sg}
\affil*{\orgdiv{Department of Industrial Systems Engineering and Management}, \orgname{National University of Singapore}}

\abstract{In this work, we present a new deterministic partition-based global optimization algorithm, HALO (Hybrid Adaptive Lipschitzian Optimization), which uses estimates of the local Lipschitz constants associated with different sub-regions of the objective function's domain to compute lower bounds and guide the search toward global minimizers. These estimates are obtained by adaptively balancing the global and local information collected from the algorithm, based on absolute slopes. HALO is hyperparameter-free, eliminating the need for manual tuning, and it highlights the most important variables to help interpret the optimization problem. We also introduce a coupling strategy with local optimization algorithms, both gradient-based and derivative-free, to accelerate convergence. We compare HALO with popular global optimization algorithms on hundreds of test functions. The numerical results are very promising and demonstrate that HALO can expand our arsenal of efficient procedures of efficient procedures for challenging real-world black-box optimization problems. The Python code of HALO is publicly available on GitHub.\footnote{\url{https://github.com/dannyzx/HALO}}}
\maketitle
\section{Introduction}
In this work, we consider the following optimization problem 
\begin{equation}
\label{eq:min_problem}
\begin{alignedat}{3}
	\min_{\mathbf{x}}  &       \,\,     f(\mathbf{x})  \\ 
	\text{subject to: }& \mathbf{x}\in \mathcal{D}   
\end{alignedat}
\end{equation}
where $f : \mathbb{R}^N \rightarrow \mathbb{R}$, $\mathbf{x}\in\mathbb{R}^N$ and $\mathcal{D}$ is the feasible region of $\mathbb{R}^N$, defined by box constraints $\mathcal{D} = \{\mathbf{x}\in \mathbb{R}^N : \mathbf{l} \leq \mathbf{x}\leq \mathbf{u} \}$. The $N$-dimensional vectors $\mathbf{l}$ and $\mathbf{u}$ represent the lower and upper bound on the decision variable $\mathbf{x}$.

Global optimization (GO) is a critical field in operations research, focusing on methodologies and algorithms for finding the global optimal solution to Problem \ref{eq:min_problem}. The pursuit of finding the global optimum covers many scientific disciplines due to its fundamental importance.
Extensive academic research has been dedicated to GO, resulting in a wide range of proposed approaches. Stochastic and evolutionary algorithms, such as controlled random search \cite{price1977controlled}, genetic algorithm \cite{davis1991handbook}, simulated annealing \cite{johnson1989optimization}, particle swarm optimization \cite{kennedy1995particle}, and covariance matrix adaptation evolutionary strategy (CMA-ES) \cite{hansen2001completely}, have proven effective in numerous real-world applications.
Another class of GO algorithms involves constructing response surfaces of the objective function using statistical models like gaussian processes \cite{jones1998efficient, mockus1978application} or radial basis functions \cite{gutmann2001radial, regis2007stochastic}. These methodologies exhibit good performance in scenarios where evaluating the objective function $f$ is computationally expensive and evaluation budget is limited. 

Lipschitz optimization, extensively studied and developed over the past decades \cite{horst2013global, horst2013handbook, sergeyev2006global, torn1989global, pinter1995global}, constitutes another class of GO algorithms.
In Lipschitz optimization, the objective function $f$ is assumed to be Lipschitz continuous on the feasible domain $\mathcal{D}$. This implies the existence of a constant $0 < L < \infty$ such that
\begin{equation}
|f(\mathbf{x}) - f(\bar{\mathbf{x}})| \leq L||\mathbf{x} - \bar{\mathbf{x}}|| \qquad \forall \mathbf{x},\bar{\mathbf{x}} \in \mathcal{D}
\end{equation}
where $L$ represents the global Lipschitz constant of the function $f$, and $||\cdot||$ denotes the Euclidean norm (alternative norms can be utilized, as explored in \cite{paulavivcius2007analysis}). 
The requirement for $f$ to be Lipschitz-continuous is not overly restrictive, as this class encompasses a wide range of functions. 
For instance, any continuously differentiable function defined within a convex and compact set is Lipschitz continuous \cite{pinter1995global}, with $L$ corresponding to the maximum norm of the gradient $L=\max_{\mathbf{x} \in \mathcal{D}} ||\nabla f(\mathbf{x})||$.
Leveraging the Lipschitz continuity of the objective function $f$ is particularly useful, as it allows for the computation of valid lower bounds across the domain $\mathcal{D}$. 
By knowing the function value at a point $f(\bar{\mathbf{x}})$, it becomes possible to compute the lower bound for any $\mathbf{x}\in\mathcal{D}$ that satisfies the following inequality
\begin{equation}
f(\bar{\mathbf{x}}) - L||\mathbf{x} - \bar{\mathbf{x}}|| \leq f(\mathbf{x})
\end{equation}
Exploiting this information algorithmically enables the design of efficient numerical procedures for seeking the global optimum solution to Problem \ref{eq:min_problem}.
\section{Related Work}\label{sec:2}
Over the years, researchers have developed various Lipschitz optimization methods (e.g., \cite{shubert1972sequential, mladineo1986algorithm}) which, given the value of the Lipschitz constant $L$, can produce valid lower bounds for the objective function $f$, guiding the numerical procedure toward the global minimizers.

The difficulty is when the value $L$ of the Lipschitz constant is not available in advance such as in the majority of real-world applications. 
To effectively overcome this issue, an estimation of it is usually carried out given the information acquired by the optimization procedure during its iterations as in \cite{sergeyev2006global} and using a diagonal partition scheme in \cite{kvasov2012lipschitz, sergeyev2017deterministic, sergeyev2015deterministic2}. 

One popular method that circumvents the need for estimating the Lipschitz constant $L$ is the DIRECT (DIviding RECTangles) algorithm \cite{jones1993-JOTA}. This derivative-free, partition-based algorithm operates under the determination of \textit{potentially optimal rectangles}, where any value of the global Lipschitz constant $L$ within the range $(0, \infty)$ is implicitly considered during the selection of the partitions to explore. DIRECT possesses several noteworthy features, including a simple partition scheme, an everywhere-dense property that ensures convergence towards the global optimal solution as the algorithm iterates indefinitely, and a minimal number of hyperparameters. In fact, DIRECT has only one hyperparameter and has exhibited robustness in practical applications \cite{jones1993-JOTA}. These characteristics have inspired the development of several similar algorithms based on the concept of potentially optimal partitions \cite{gablonsky2001, mockus_2011, paulavivcius2020globally, mockus2017application, liu2015mrdirect, liu2017improving, lera2015deterministic, paulavivcius2014simplicial}.

Other noteworthy algorithms in the Lipschitz optimization class focus on estimating the local Lipschitz constant $L_i$ for each sub-region $\mathcal{D}_i$ of the domain $\mathcal{D}$. In \cite{sergeyev1995information}, an algorithm is proposed for optimizing univariate functions, while for multivariate functions, a dimensionality reduction method employing Peano curves is utilized \cite{strongin1992algorithms, sergeyev2013introduction}. The estimation of local Lipschitz constants is employed in a diagonal partition-based algorithm described in \cite{kvasov2003local}. 
In \cite{gergel2015local}, the computation of local Lipschitz constants is integrated into a scheme that reduces the initial multidimensional problem to a set of recursively and adaptively connected univariate subproblems, as further explored in \cite{gergel2016adaptive}. A similar concept is proposed when the first derivative is available in \cite{gergel1997global}. However, it is worth noting that the performance of these approaches is highly susceptible to a reliability parameter that users must define beforehand, as emphasized in \cite{kvasov2003local, gergel2015local, jones2020direct, sergeyev2015deterministic}. 
At the same time, even with a misspecified reliability parameter, these methods can exhibit stronger theoretical convergence guarantees than DIRECT-type methods \cite{Sergeyev01011998, sergeyev2017deterministic}.
Similar approaches have been developed specifically for univariate objective functions \cite{sergeyev1998global, sergeyev2021novel, posypkin2022efficient, kvasov2009univariate}.

Furthermore, successful hybridizations have been achieved by combining some of the aforementioned approaches with local optimization methods. For instance, the DIRECT algorithm has been hybridized with a truncated Newton method in \cite{liuzzi2010direct} and a coordinate descent approach in \cite{liuzzi2016exploiting}, resulting in improved performance, as also observed in the methodology presented in \cite{paulavivcius2020globally}.

\section{Main Contribution}
In this work, we present a novel GO algorithm that incorporates an adaptive procedure for estimating the local Lipschitz constant associated with each partition $\mathcal{D}_i$. 
Our estimation is obtained through a convex combination of the absolute variation of the objective function around the partition $\mathcal{D}_i$ and the current estimate of the global Lipschitz constant $L$. The coefficients in the convex combination are dynamically computed based on the size of the partition $\mathcal{D}_i$, eliminating the need for users to define any crucial hyperparameters in advance.

We propose a criterion for selecting the most promising partitions based on the information provided by the local Lipschitz constants. This criterion deviates from the concept of potentially optimal hyperrectangles employed in the DIRECT algorithm.
Furthermore, we establish the asymptotic convergence of our algorithm called HALO (Hybrid Adaptive Lipschitzian Optimization) to a global minimum, ensuring the everywhere dense property. 
To expedite the convergence speed of HALO towards stationary points, we introduce a simple coupling strategy with two different local optimization algorithms: a gradient-based and a derivative-free optimization approach presented respectively in \cite{byrd1995limited} and \cite{lucidi2002-COA}.

Beyond its optimization capabilities, HALO provides the opportunity to gain valuable insights from the black-box objective function. Upon terminating the optimization process, we show that HALO allows the extraction of information about the variables that most significantly influenced the variation of the objective function. This feature is particularly important in scenarios where problem interpretation is necessary.

We evaluate the performance of HALO by comparing it to popular GO algorithms such as DIRECT \cite{jones1993-JOTA}, CMA-ES \cite{hansen2001completely}, L-SHADE \cite{tanabe2014improving}, PLOR \cite{mockus2017application} and the hybridization of DIRECT known as DIRMIN \cite{liuzzi2016exploiting}. To assess the algorithms, we conduct an extensive numerical campaign using well-known test functions collected from \cite{Jamil_2013}, as well as two additional test function generators specifically designed for GO problems presented in \cite{schoen1993wide, gaviano2003algorithm}. In total, we consider over one thousand test functions.

The numerical results demonstrate that HALO exhibits strong competitiveness and significantly enhances our arsenal of efficient procedures for tackling challenging real-world GO problems.
\section{HALO: Hybrid Adaptive Lipschitzian Optimization}
\begin{algorithm}
  \DontPrintSemicolon
\Function{main(\text{stopping\_criteria}, $\mathcal{D}$)}{
    $\mathcal{I}_0=\emptyset$;
    $\mathcal{L}_0 = \emptyset$\;
    \For{$k$ \KwTo \text{max\_iter}}{
    $\mathcal{I}^{\star}_{k} \assign$ \FuncCall{selection}{$\mathcal{I}_{k}$, $\mathcal{L}_{k}$} \;
    \For{$i^{\star}_k$ \KwIn $\mathcal{I}^{\star}_{k}$}{
        \For{$j$ \KwTo $J$}{
            $\mathbf{x}^{j}_{i^{\star}_{k}} \assign \mathbf{x}^{j}_{i^{\star}_{k}} \in \mathcal{D}_{i^{\star}_{k}}$\;
            $f(\mathbf{x}^{j}_{i^{\star}_{k}}) \assign $ \FuncCall{fun\_eval}{$\mathbf{x}^{j}_{i^{\star}_{k}}$}\;
            }
        $\mathcal{I}_{k} \assign$ \FuncCall{partitioning}{$\mathcal{D}_{i^{\star}_{k}}$} \;
        $\mathcal{L}_{k} \assign$ \FuncCall{compute\_local\_Lipschitz}{$\mathcal{I}_{k}$, $\mathcal{D}_{k}$} \;
    }
    }
\Return $\min_{i_k \in \mathcal{I}_k} f(\mathbf{x}_{i_k})$\;
}
\Function{selection($\mathcal{I}_{k}$)}{
        $\mathcal{I}^{\star}_{k} \assign \mathcal{I}^{\star}_{k} \subseteq \mathcal{I}_{k}$\;
        \Return $\mathcal{I}^{\star}_{k}$\;
  }
\Function{partitioning($\mathcal{D}_{i^{\star}_{k}}$)}{
$\mathcal{D}_{i_k^\star} \assign \bigcup\limits_{j =  1}^{J} \mathcal{D}_{i_k^\star}^{j}$ \;
$\mathcal{I}_{k} \assign \mathcal{I}_k \cup \{|\mathcal{I}_k| + j\}_{j=1}^{J}$ \;
\Return $\mathcal{I}_{k}$\;
}
\caption{General Partition Algorithm with Local Lipschitz Constants Estimates}
\label{alg:general}
\end{algorithm}
Before delving into the detailed workings of HALO, it is necessary to provide an overview of its main steps and introduce some notation. Algorithm \ref{alg:general} provides a high-level description of the HALO algorithm, outlining its main components and the flow of operations as it consists of three parts: a main program, a selection step, and a partitioning step.

The \texttt{Main} function, which takes the maximum number of iterations (\texttt{MaxIter}) and the box domain ($\mathcal{D}$) as input, serves as the main entry point of the algorithm. 
In line 3, the algorithm begins with a loop over the iterations $k$. Within this loop, the \texttt{Selection} function is called at line 4 to select partitions based on certain criteria, which make use of the Lipschitz constant information of each partition (not yet specified).

The selected partitions are represented by the subset of indices $\mathcal{I}^{\star}_{k}$ obtained from the set of indices $\mathcal{I}_{k}$, where each index $i_k$ corresponds to a partition $\mathcal{D}_{i_k}$ at iteration $k$. 
Lines 5 and 6 contain two nested 'for' loops. 
The outer loop iterates over each index $i^{\star}_k$ of the selected partitions, while the inner loop (line 6) performs $J$ function evaluations within each partition.
The function evaluations are carried out at points $\mathbf{x}^{j}_{i^{\star}_{k}}$ sampled from the partition $\mathcal{D}_{i^{\star}_{k}}$. 
The corresponding function values $f(\mathbf{x}^{j}_{i^{\star}_{k}})$ are computed using the \texttt{Fun\_Eval}. 
Once the inner loop completes, the set $\mathcal{D}_{i^{\star}_{k}}$ is divided into $J$ non-overlapping subsets (line 16), denoted as $\mathcal{D}_{i_k^\star}^{j}$ for $j = 1, \dots, J$. Furthermore, the local Lipschitz constant around every partition is estimated at line 10 and used in the \texttt{Selection} function. 
In line 17, the set of indices $\mathcal{I}_{k}$ is updated to include the new indices resulting from the partitioning step. \\
The entire process described above is repeated until the termination criteria are met.

\subsection{Selection}
In the previous section, we presented a general partition scheme for Lipschitz optimization. Now, we dig into the algorithm's details by introducing the steps for estimating the local Lipschitz constants and selecting the partitions based on this information.

\subsubsection{Adaptive Estimation of the Local Lipschitz Constants}
In Lipschitz optimization, it is common practice to obtain information about the Lipschitz constant $L$ and compute lower bounds of the objective function $f$ to guide the algorithm towards regions where these bounds are lower. However, relying solely on the global Lipschitz constant estimates may not accurately reflect the behavior of the objective function in specific regions of the domain $\mathcal{D}$. Beyond global estimates, it is well established that employing local Lipschitz constants $\tilde{L}_{i}$ for each subregion $\mathcal{D}_i$ can improve efficiency; see Section \ref{sec:2}, for example, the local tuning approaches developed in \cite{sergeyev1998global, sergeyev2017deterministic, sergeyev1995information, strongin2013global, gergel2015local}.

To better understand why local estimates can be more informative, consider a scenario where the objective function $f$ is mostly flat, except for a small region within $\mathcal{D}$ where it exhibits highly chaotic behavior (i.e., a much higher Lipschitz constant). In such cases, using a global Lipschitz constant to compute lower bounds across all of $\mathcal{D}$ would lead to poor estimates in the flatter regions. This is because the global constant is dominated by the steep variations in that small chaotic area, overestimating the possible change of $f$ elsewhere and thus weakening the bounds.
Therefore, it is more informative to estimate the local Lipschitz constant $\tilde{L}_{i}$ associated with each partition $\mathcal{D}_i$ instead of relying entirely on the global Lipschitz constant $L$.

To introduce a local source of information about the function $f$ when estimating the lower bounds, we consider the norm of the gradient at the centroid $\mathbf{x}_{i}$. However, to maintain the derivative-free nature of our procedure and to not consume too many function evaluations, we approximate the gradients around each centroid $\mathbf{x}_{i}$ of the partitions $\mathcal{D}_i$ using the points sampled thus far by the algorithm.
This algorithmic choice has the effect that the accuracy of these gradient approximations directly depends on the size of each partition.

As the size of a partition $\mathcal{D}_{i}$ decreases, the uncertainty about the behavior of the objective function around $\mathcal{D}_{i}$ decreases as well. Therefore, it is reasonable to balance our estimate of the local Lipschitz constant $\tilde{L}_{i}$ based on the size of its corresponding partition $\mathcal{D}_{i}$. To achieve this, we propose an adaptive formula

\begin{definition}
Given a subregion $\mathcal{D}_{i_k}$ of the unit hypercube domain $\mathcal{D}$ at iteration $k$, let
\[
\alpha_{i_k} = \frac{||\mathbf{u}_{i_k} - \mathbf{l}_{i_k}||}{\sqrt{N}}, \quad \alpha_{i_k} \in (0, 1)
\]
and define the global Lipschitz constant estimate as
\[
\tilde{L}_k = \max_{i_k \in \mathcal{I}_k} ||\widetilde{\nabla} f(\mathbf{x}_{i_k})||
\]
Then we define
\begin{equation}\label{eq:formula_adaptive}
\tilde{L}_{i_k} = \alpha_{i_k} \tilde{L}_k + (1 - \alpha_{i_k}) ||\widetilde{\nabla} f(\mathbf{x}_{i_k})||
\end{equation}
\end{definition}

In Eq. \ref{eq:formula_adaptive}, the term $\alpha_{i_k}$ is the ratio between the diagonal of the partition $\mathcal{D}_{i_k}$ and the main diagonal of the feasible set $\mathcal{D}$ (normalized to a unit hypercube). As the diagonal of the partition $\mathcal{D}_{i_k}$, denoted by $||\mathbf{u}_{i_k} - \mathbf{l}_{i_k}||$, decreases, the local information in the neighborhood of $\mathcal{D}_{i_k}$ becomes more precise, and the estimate of the local Lipschitz constant is primarily influenced by the norm of the gradient approximation at $\mathbf{x}_{i_k}$, namely $||\widetilde{\nabla}{f(\mathbf{x}_{i_k})}||$. 
On the other hand, if the diagonal of the partition $\mathcal{D}_{i_k}$ is approximately equal to the maximum diagonal, the local information may not be reliable, and the estimate of the local Lipschitz constant relies more on the estimate of the global Lipschitz constant $\tilde{L}_{k}$, which represents the maximum absolute variation observed during the iterations of the algorithm.
In this way, our formula (Eq. \ref{eq:formula_adaptive}) naturally balances global and local information in a self-adaptive manner, eliminating the need for user-defined parameters as traditional local tuning techniques \cite{sergeyev1998global, sergeyev2017deterministic, sergeyev1995information, strongin2013global, gergel2015local}.

In the following section, we define the criteria for selecting the partitions $\mathcal{D}_{i_k}$ based on the information collected from our estimates of the local Lipschitz constants.
\subsubsection{Selection Criteria}\label{subsubsection:select_rules}
\begin{figure}[htb!]
	\centering
	\begin{subfigure}[t]{0.52\textwidth}
	\includegraphics[width=\textwidth]{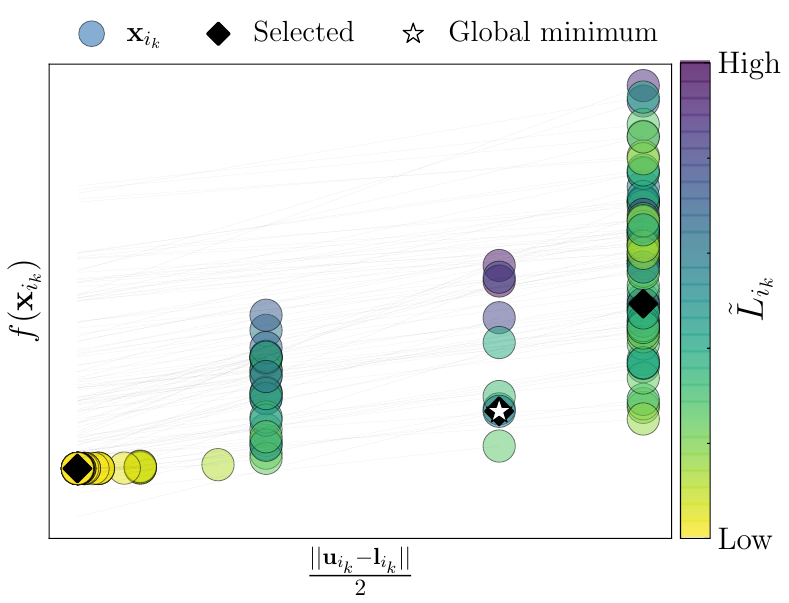}
	\caption{HALO.}
        \label{fig:selection_HALO}
	\end{subfigure}
         \begin{subfigure}[t]{0.45\textwidth}
	\includegraphics[width=\textwidth]{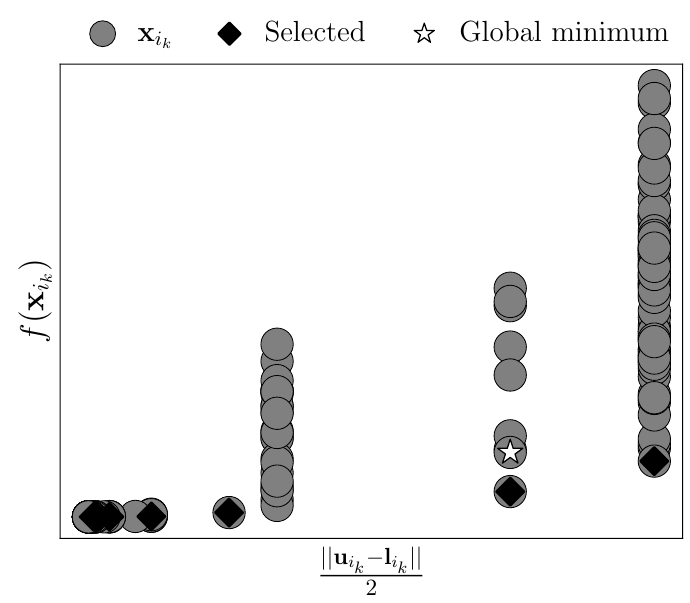}
	\caption{DIRECT.}
         \label{fig:selection_DIRECT}
	\end{subfigure}
	\caption{Comparison between the selection criteria of HALO (Fig. \ref{fig:selection_HALO}) and the potentially optimal selection criteria of DIRECT (Fig. \ref{fig:selection_DIRECT}). The $x$-axis represents the distance from the center to the vertices of each partition, while the $y$-axis represents the corresponding function values. In Fig. \ref{fig:selection_HALO}, the colormap indicates the values of the local Lipschitz constants $\tilde{L}_{i_k}$. The gray lines in Fig. \ref{fig:selection_HALO} emphasize the lower bound values obtained at the vertices, with their slope determined by $\tilde{L}_{i_k}$. The selected partitions and the location of the global minimum are denoted by black rhombuses and a white star respectively.}
\label{fig:selection}
\end{figure}

In the following definition, we show which partitions $\mathcal{D}_{i_k}$ are selected in HALO.
\begin{definition}\label{def:selection}
	Given the estimation of the local Lipschitz constant $\tilde{L}_{i_k}$ for every $i_k \in \mathcal{I}_{k}$ from Eq. \ref{eq:formula_adaptive}, then the partition $\mathcal{D}_{i_k}$ is selected if at least one of the following conditions is satisfied
	\begin{enumerate}[align=left,labelindent=11cm,labelwidth=1.5cm,label={Criterion ~\arabic*},itemindent=0em,leftmargin=2cm]
		\item \label{choice 2} Given $i_k \in \mathcal{I}_k$, then 
		\begin{equation}
			f(\mathbf{x}_{i_k}) - \tilde{L}_{i_k} \frac{||\mathbf{u}_{i_k} - \mathbf{l}_{i_k}||}{2} \leq f(\mathbf{x}_{j_k}) - \tilde{L}_{j_k} \frac{||\mathbf{u}_{j_k} - \mathbf{l}_{j_k}||}{2} \qquad \forall j\in \mathcal{I}_{k} \label{eq:c21} \\
		\end{equation}
		\item \label{choice 3} Given $i_k \in \mathcal{I}_k$, then
		\begin{equation}
			f(\mathbf{x}_{i_k}) = f_{\min} \label{eq:c32}
		\end{equation}
	where $f_{\min} = \min_{j_k \in \mathcal{I}_k} f(\mathbf{x}_{j_k})$
		\item \label{choice 1} Given $i_k \in \mathcal{I}_k^{\text{max}}$ and 
		\begin{equation}
			f(\mathbf{x}_{i_k}) - \tilde{L}_{i_k} \frac{||\mathbf{u}_{i_k} - \mathbf{l}_{i_k}||}{2} \leq f(\mathbf{x}_{j_k}) - \tilde{L}_{j_k}\frac{||\mathbf{u}_{j_k} - \mathbf{l}_{j_k}||}{2} \qquad \forall j\in \mathcal{I}_k^{\text{max}} \label{eq:c11}
		\end{equation}
		where $\mathcal{I}_k^{\text{max}} = \{i_k \in \mathcal{I}_{k} : \frac{||\mathbf{u}_{i_k} - \mathbf{l}_{i_k}||}{2} = \max_{i_k \in \mathcal{I}_k}\frac{||\mathbf{u}_{i_k} - \mathbf{l}_{i_k}||}{2}\}$.
	\end{enumerate}
\end{definition}
The selection criteria defined in Def. \ref{def:selection} determine which partitions $\mathcal{D}_{i_k}$ are selected in the HALO algorithm. Let's examine each criterion
\begin{itemize}
\item \ref{choice 2}: This criterion selects the partition $\mathcal{D}_{i_k}$ that achieves the lowest lower bound, considering the information collected from the local Lipschitz constants. The lower bound is computed as $f(\mathbf{x}_{i_k}) - \tilde{L}_{i_k} ||\mathbf{u}_{i_k} - \mathbf{l}_{i_k}||/2$, where $f(\mathbf{x}_{i_k})$ is the function value at the centroid of $\mathcal{D}_{i_k}$, $\tilde{L}_{i_k}$ is the estimated local Lipschitz constant, and $||\mathbf{u}_{i_k} - \mathbf{l}_{i_k}||/2$ is the distance from the center to the vertices of the partition. This criterion aims to select partitions that provide the most promising lower bounds based on the local behavior of the objective function.

\item  \ref{choice 3}: This criterion ensures that the partition with the lowest function value, denoted as $f_{\min}$, is always selected. This condition is beneficial when HALO is used in high-dimensional spaces and in a limited number of function evaluation scenarios, allowing the algorithm to obtain a decent reduction in the objective function value in a few iterations.

\item  \ref{choice 1}: Among the partitions with the maximum diagonal size, which are denoted as $\mathcal{I}_k^{\text{max}}$, this criterion selects the partition that achieves the lowest lower bound. The lower bound is computed similarly as in  \ref{choice 2}. This criterion is significant for the convergence of the algorithm, as we will see in the following sections.
\end{itemize}
In Fig. \ref{fig:selection}, we can compare the selection criteria of HALO with those of the DIRECT algorithm, specifically targeting the potentially optimal partitions. 
We used the synthetic function generator presented in \cite{schoen1993wide} for this particular example.
Notably, HALO follows distinct trajectories from DIRECT, as none of its selected partitions align with the potentially optimal rectangles. 
In fact, HALO has the ability to choose rectangles that can be arbitrarily distant from the potentially optimal partitions of DIRECT. This allows HALO to proactively explore interesting partitions by leveraging its knowledge of the local Lipschitz constants. As illustrated in Fig. \ref{fig:selection}, HALO efficiently identifies the partition (indicated by the white star) where the global minimum is located since is the one that satisfies \ref{choice 2}, while DIRECT experiences a delay in identifying the optimal partition.
\subsection{Partitioning}\label{subsection:Sampling}
In this section, we discuss the geometric properties of the feasible set $\mathcal{D}$ and its partitioning, as well as the adaptive approximation of derivatives around each partition during the iterations.
\subsubsection{Division and Sampling}\label{subsubsection:DIRECT}
During the global search carried out by HALO, we employ the same sampling and division criteria as the DIRECT algorithm \cite{jones1993-JOTA}. This choice was made due to the robustness and relative simplicity of the DIRECT algorithm, making it easier to understand and implement the code.

Let's describe how the partitions $\mathcal{D}_{i_{k}}$ are generated in the DIRECT algorithm. We define the vector $\mathbf{s}_{i_{k}} \in \mathbb{R}^N$ associated with the partition $\mathcal{D}_{i_{k}}$ which measures the distance from its centroid to its boundary along each coordinate direction. Specifically, $\mathbf{s}_{i_{k}}$ represents half the length of the sides of the partition at iteration $k$
\begin{equation}\label{eq:sides}
\mathbf{s}_{i_{k}} = [|u_{i_k}^n - l_{i_k}^n|/2]_{n=0}^{N}    
\end{equation}

we define the set of indices $\mathcal{P}_{i_{k}}$ that contains the $n$th coordinate directions where the partition $\mathcal{D}_{i_{k}}$ has its longest side
\begin{equation}\label{eq:set_P}
\mathcal{P}_{i_{k}} = \{n \in \mathbb{N} : s_{i_{k}}^{max} =\max\{s_{i_{k}}^n\}_{n=0}^{N}\}
\end{equation} 
Subsequently, two points $\mathbf{x}_{i_{k}}^{p_1}$ and $\mathbf{x}_{i_{k}}^{p_2}$ are sampled for each coordinate axis parallel to the longest side of $\mathcal{D}_{i_{k}}$ (namely along the $p$th coordinate)
\begin{equation}\label{eq:sampling}
	\mathbf{x}_{i_{k}}^{p_1} = \mathbf{x}_{i_{k}} + \Delta_{i_{k}}\mathbf{e}_p, \qquad \mathbf{x}_{i_{k}}^{p_2} = \mathbf{x}_{i_{k}} - \Delta_{i_{k}}\mathbf{e}_p, \qquad \forall p \in \mathcal{P}_{i_{k}}
\end{equation} 
here, $\mathbf{e}_p$ represents the $p$th orthonormal base of $\mathbb{R}^N$ and the scalar $\Delta_{i_{k}}$ is given by two third the longest side of $\mathcal{D}_{i_{k}}$
\begin{equation}\label{eq:delta_sampling}
	\Delta_{i_{k}} = \frac{2}{3} s_{i_{k}}^{max}
\end{equation} 
From Eq. \ref{eq:sampling} is evident that two new points are generated along the coordinate of the longest side of the selected partition, resulting in a total of $J = 2 |\mathcal{P}_{i_{k}}|$ new points generated per selected partition.

Once the sampling within a selected partition $\mathcal{D}_{i_{k}}$ is complete, we proceed with the division. In the DIRECT algorithm, the idea is to divide a partition into hypercubes and hyperrectangles, ensuring that the partition with the lowest objective function value always has the longest diagonal. \\
We define the set $\mathcal{T}_{i_{k}}$ which collects the minimum objective function value among the two newly generated points $\mathbf{x}_{i_{k}}^{p_1}$ and $\mathbf{x}_{i_{k}}^{p_2}$ for each coordinate axis parallel to the longest side of $\mathcal{D}_{i_{k}}$
\begin{equation}\label{eq:set_T}
	\mathcal{T}_{i_{k}} = \mathcal{T}_{i_{k}} \cup \{\min \{f(\mathbf{x}_{i_{k}}^{p_1}), f(\mathbf{x}_{i_{k}}^{p_2}\}\}, \qquad \forall p \in \mathcal{P}_{i_{k}}
\end{equation}
we proceed to divide $\mathcal{D}_{i_{k}}$ in the order given by $\mathcal{T}_{i_{k}}$ and perpendicular to the direction $\mathbf{e}_p$ such that all the $p$th sides of $\mathbf{s}_{i_{k}}$ are reduced by half of $\Delta_{i_{k}}$ (or equivalently by a third the longest side of the partition $\mathcal{D}_{i_{k}}$)
\begin{equation}
    s_{i_{k}}^{p} = \frac{1}{2} \Delta_{i_{k}} = \frac{1}{3} s_{i_{k}}^{max},  \qquad s_{i_{k}}^{p_1} = s_{i_{k}}^{p_2} = s_{i_{k}}^{p} \qquad \forall p \in \mathcal{P}_{i_{k}}
\end{equation}
Additionally, the sides $s_{i_{k}}^{p_1}$ and $s_{i_{k}}^{p_2}$ of the two new partitions are updated accordingly.
\begin{figure}[htb!]
	\centering
	\begin{subfigure}[t]{0.32\textwidth}
		\includegraphics[width=\textwidth]{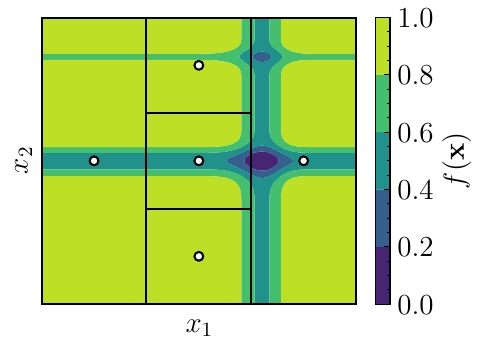}
		\caption{$k=0$}
		\label{fig:part_k0}
	\end{subfigure}
	\begin{subfigure}[t]{0.32\textwidth}
		\includegraphics[width=\textwidth]{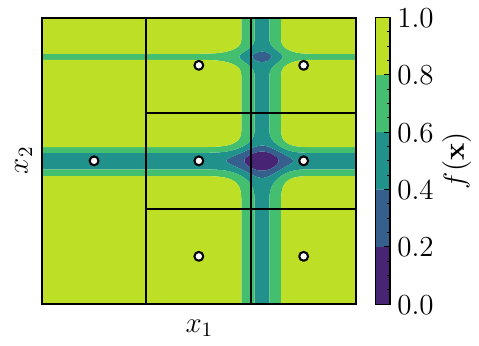}
		\caption{$k=1$}
	    \label{fig:part_k1}
	\end{subfigure}
	\begin{subfigure}[t]{0.32\textwidth}
		\includegraphics[width=\textwidth]{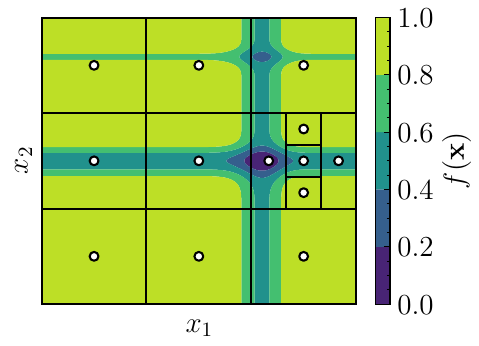}
		\caption{$k=2$}
		\label{fig:part_k2}
	\end{subfigure}
	\caption{A graphical representation of the division and sampling step. White dotted points the centroids $\mathbf{x}_{i_{k}}$, in black we highlight the boundaries of the partitions. The colorbar shows normalized function values.}
	\label{fig:halo_part}
\end{figure}
\subsubsection{Adaptive Gradients Approximation}\label{subsubsection:adagrad}
In this section we describe how the vector $\widetilde{\nabla}{f(\mathbf{x}_{i_k})}$ from Eq. \ref{eq:formula_adaptive} is computed.
After the sampling and division processes, the algorithm collects information about the behavior of the objective function around each partition. 
To achieve this, a vector is associated with each partition, which describes the intensity of variation of the objective function inside that partition in terms of absolute slopes. 
This vector adaptively collects the absolute slopes along each coordinate computed based on the points generated so far by the algorithm. 
We define the nearest neighbors of a partition $\mathcal{D}_{i_k}$ as the indices of the nearest points along each coordinate axis.

%
\begin{definition}\label{def:nearest_neighbors}
Given the partition $\mathcal{D}_{i_k}$ and its centroid $\mathbf{x}_{i_k}$ we define the nearest neighbors of $\mathbf{x}_{i_k}$ the set of indices $\mathcal{N}_{i_k}$ which contains the indices of the nearest points along the coordinate axis of $\mathbb{R}^N$.
\end{definition}
The following property is satisfied by the current partition scheme
\begin{property}\label{property:dist1}
Given the representative point $\mathbf{x}_{i_k}$ which is center of its relative partition $\mathcal{D}_{i_k}$. Then if for every iteration $ k>\bar{k}$, the partition $\mathcal{D}_{i_k}$ is selected, the nearest neighbours of $\mathbf{x}_{i_k}$ are given by
    \begin{equation}
    	\mathbf{x}_{i_k} \pm \Delta_{i_k}\mathbf{e}_n \in \mathcal{N}_{i_k}, \qquad  n=1, \dots ,N
    \end{equation}
	where $\mathbf{e}_n$ is the $n$th orthonormal basis of $\mathbb{R}^N$ and $\mathcal{N}_{i_k}$ the set of the nearest neighbours of $\mathbf{x}_{i_k}$.
\end{property}
The nearest neighbors of a partition $\mathcal{D}_{i_k}$ are defined as the indices of the nearest points along each coordinate axis. If the partition $\mathcal{D}_{i_k}$ is always selected after a certain iteration $\bar{k}$, the nearest neighbors of its centroid $\mathbf{x}_{i_k}$ are located at distances $\pm \Delta_{i_k}$ along each coordinate axis.


If a selected partition $\mathcal{D}_{i_k}$ is a hypercube, it generates $J = 2N$ nearest points along all $N$ coordinate axes. In this case, we compute the absolute slope at the centroid $\mathbf{x}_{i_k}$ using a central difference formula
\begin{equation}\label{eq:central_diff}
	|\widetilde{\nabla}{f(\mathbf{x}_{i_k})}_n| = \frac{|f(\mathbf{x}_{i_k}^{p_1}) - f(\mathbf{x}_{i_k}^{p_2})|}{2\Delta_{i_k}}, \qquad  n = 1, \dots, N
\end{equation}
Here, $\mathbf{x}_{i_k}^{p_1}$ and $\mathbf{x}_{i_k}^{p_2}$ are the two points generated from the centroid $\mathbf{x}_{i_k}$ along the $n$th coordinate axis.
Similarly, we need to compute the slopes for the new points $\mathbf{x}_{i_{k}}^j$ generated by the partition $\mathcal{D}_{i_k}$ at iteration $k$. The approximation of the gradient associated with $\mathcal{D}_{i_k}^j$ is given by

\begin{equation}
	|\widetilde{\nabla}{f(\mathbf{x}_{i_k}^j)}_n| = \begin{cases}
		\frac{|f(\mathbf{x}_{i_k}^j) - f(\mathbf{x}_{i_k})|}{\Delta_{i_k}},  & \text{if $\mathbf{x}_{i_k} \in \mathcal{N}_{i_k}^j$} .\\
		|\widetilde{\nabla}{f(\mathbf{x}_{i_k})}_n|,  & \text{otherwise}.
	\end{cases}
\end{equation}
For the remaining $N-1$ coordinates, the absolute slopes $|\widetilde{\nabla}{f(\mathbf{x}_{i_k}^j)}_n|$ are set to the slope of $\mathbf{x}_{i_{k}}$. Although this approach may lead to an incorrect estimation of the gradient associated with the partition $\mathcal{D}_{i_{k}}^j$, the effect is negligible as either the partition $\mathcal{D}_{i_{k}}^j$ will be selected and its gradient updated or the size of the partition $\mathcal{D}_{i_{k}}^j$ is already relatively small.

If the selected partition $\mathcal{D}_{i_k}$ is a hyperrectangle, it generates $J = 2|\mathcal{P}_{i_k}|$ nearest points along all $|\mathcal{P}_{i_k}|$ coordinate axes. The absolute slope is calculated using the same central difference formula as before
\begin{equation}
	|\widetilde{\nabla}{f(\mathbf{x}_{i_k})}_n| = \begin{cases}
		\frac{|f(\mathbf{x}_{i_k}^{p_1}) - f(\mathbf{x}_{i_k}^{p_2})|}{2\Delta_{i_k}},  & \text{if $\mathbf{x}_{i_k}^{p_1}, \mathbf{x}_{i_k}^{p_2} \in \mathcal{N}_{i_k}$} .\\
		|\widetilde{\nabla}{f(\mathbf{x}_{i_k})}_n|,  & \text{otherwise}.
	\end{cases}
\end{equation}
For the remaining $N-|\mathcal{P}_{i_{k}}|$ coordinates we leave the slope unchanged in those directions. 
As shown previously we also update the slopes for the points that $\mathcal{D}_{i_k}$ has generated. In this case, the procedure is the same as in the case if the partition $\mathcal{D}_{i_k}$ is a hypercube when selected.

The simple procedure described above tries to adaptively update the absolute slopes at each selected centroid and at its neighbors only using the points sampled so far by the algorithm. 
\subsection{Coupling Strategy with Local Optimization Methods}\label{section:coupling strategy}
In this work, we propose to integrate our approach with local optimizers to improve its convergence speed. The reason behind this is that the sampling scheme used in the DIRECT algorithm, as described in section \ref{subsubsection:DIRECT}, can be inefficient when the algorithm approaches a stationary point. The rigid structure of the sampling scheme, where the algorithm is forced to sample along predefined coordinate directions with constant step sizes, can slow down the convergence of the HALO algorithm.

To address this issue, we introduce a coupling strategy that incorporates local optimization routines only when a point is likely to be in the vicinity of a stationary point, to avoid overusing them. The strategy is as follows: if the selected partition satisfies either \ref{choice 2} or \ref{choice 3} (or both) in Def. \ref{def:selection}, and if the half diagonal of the partition is less than a scalar $\beta$ then the centroid of the partition is considered as a starting point for the local optimizer

%
\begin{equation}\label{eq:local_beta}
	\frac{||\mathbf{u}_{i_k} - \mathbf{l}_{i_k}||}{2}\leq \beta
\end{equation}
%
To ensure that the local optimizer doesn't start from a similar location in the domain, we perform a nearest-neighbor search. Before starting the local optimization from the centroid, we identify all the points within a ball centered at the centroid with radius $r$. These points are saved in a set called $\mathcal{C}_k$, which includes the indices $j$ of the neighboring points. The set $\mathcal{C}_k$ is defined as follows

%
\begin{equation}
	\mathcal{C}_k = \{ j_k : ||\mathbf{x}_{i_k} - \mathbf{x}_{j_k}|| \leq r\}, \qquad \forall j_k \in I_k, \qquad i_k \neq j_k
\end{equation}
Then, the local optimizer starts from the centroid $\mathbf{x}_{i_k}$ with the updated set of indices
\begin{equation}
\mathcal{C}_k \cup i_k    
\end{equation}
If a new candidate point $\mathbf{x}_{p_k}$ satisfies either condition \ref{choice 2} or \ref{choice 3} (or both), and the corresponding partition $\mathcal{D}_{p_k}$ satisfies Eq. \ref{eq:local_beta}, the local search starts from $\mathbf{x}_{p_k}$ only if the set $\mathcal{B}_k$ is empty. The set $\mathcal{B}_k$ is defined as
\begin{equation}
	\mathcal{B}_k = \{ j_k : ||\mathbf{x}_{p_k} - \mathbf{x}_{j_k}|| \leq r\}, \qquad \forall j_k \in C_k
\end{equation}
If $\mathcal{B}_k$ is empty, indicating that $\mathbf{x}_{p_k}$ is not close to any neighbors already collected in $\mathcal{C}_k$, the local optimizer starts from $\mathbf{x}_{p_k}$ and the set $\mathcal{C}_k$ is updated with the neighbors of $p_k$
\begin{equation}
	\mathcal{C}_k = \mathcal{C}_k \cup \{ j_k : ||\mathbf{x}_{p_k} - \mathbf{x}_{j_k}|| \leq r\}, \qquad \forall j_k \in I_k, \qquad p_k \neq j_k
\end{equation}
Then, the local optimizer is initiated from $\mathbf{x}_{p_k}$ and the set $\mathcal{C}_k$ is updated with $p_k$. On the other hand, if $\mathcal{B}_k$ is not empty, indicating that $\mathbf{x}_{p_k}$ is close to some neighbors in $\mathcal{C}_k$, the local search does not start from $\mathbf{x}_{p_k}$. The set $\mathcal{C}_k$ is still updated with $p_k$
\begin{equation}
\mathcal{C}_k \cup p_k    
\end{equation} 
In this case, the partition $\mathcal{D}_{p_k}$ is not sampled or divided, and the radius $r$ should be set to a relatively small value, such as $r = 10^{-4}$, to avoid excluding large portions of the domain.
\subsection{Convergence and Stopping Criteria}
In this paragraph, we will discuss the convergence properties of the HALO algorithm.
Like other DIRECT type approaches, HALO belongs to the class of adaptive partition \cite{pinter1992convergence} and divide the best methods \cite{Sergeyev01011998}. These algorithms typically ensure an everywhere dense type of convergence \cite{torn1989global}, which means they can converge to any point in the feasible region, as long as the objective function is continuous. Some algorithms in this class (see Section \ref{sec:2}) rely on a reliability parameter for estimating Lipschitz constants. Although a poorly chosen value can impact performance, these methods still enjoy strong theoretical convergence guarantees, even under imperfect parameter tuning \cite{Sergeyev01011998, sergeyev2017deterministic}.
HALO avoids the need for this parameter altogether by dynamically estimating local Lipschitz constants at each iteration, together with selecting the largest hyperrectangle, eliminating the need for the user to choose a hyperparameter that affects convergence. This design choice prevents the algorithm from getting stuck in local minima and ensures everywhere dense convergence.

We introduce a property that characterizes the partitions generated during the iterations
\begin{property}\label{property:nested}
     The sequence of partitions $\{\mathcal{D}_{i_k}\}$ is strictly nested if and only if 
	\begin{equation}
		\bigcap\limits_{k=0}^{+\infty} \mathcal{D}_{i_k} = \mathbf{x}_{i_k}
	\end{equation}
	or equivalently that
	\begin{equation}
		\lim_{k\to+\infty} ||\mathbf{u}_{i_k} - \mathbf{l}_{i_k}|| = 0
	\end{equation}
\end{property}
Property \ref{property:nested} states that when a partition $\mathcal{D}_{i_k}$ is continuously selected and divided, it collapses to $\mathbf{x}_{i_k}$. 

We can establish a proposition that shows how the everywhere dense property of HALO is achieved when the largest partitions are always selected
\begin{proposition}\label{prop:every_dense}
	Given that Property \ref{property:nested} is satisfied, if for every iteration $k$ the largest partitions $\mathcal{D}_{i_k}$ is always selected such that $i_k \in  \mathcal{I}_k^{\text{max}}$, then for $k\to \infty$ and for every $\tilde{\mathbf{x}} \in \mathcal{D}$, the sequence of sets $\{\mathcal{D}_{i_k}\}$ are dense such that
		\begin{equation}
			\bigcap\limits_{k=0}^{\infty} \mathcal{D}_{i_{k}} = \{\tilde{\mathbf{x}}\} \qquad \forall i_k \in  \mathcal{I}_k 
		\end{equation}
\begin{proof}
    The proof follows from the fact that the following identity holds for every $i_k$ and that \ref{choice 1} is always met
\begin{equation}
    \lim_{k\to+\infty} ||\mathbf{u}_{i_k} - \mathbf{l}_{i_k}|| \leq ||\mathbf{u}_{j_k} - \mathbf{l}_{j_k}||= 0
\end{equation}
with $j_k \in \mathcal{I}_k^{\text{max}}$, with $\mathcal{I}_k^{\text{max}} = \{i_k \in I_k : ||\mathbf{u}_{i_k} - \mathbf{l}_{i_k}|| = \max_{i_k \in \mathcal{I}_k}||\mathbf{u}_{i_k} - \mathbf{l}_{i_k}||\}$.
\end{proof}
\end{proposition}
Therefore, if the algorithm produces dense sets across all the domain $\mathcal{D}$, then for $k$ sufficiently large the algorithm will converge also at the global minimum of $f$. 

Generally, rather than producing an everywhere dense set we often aim for an algorithm that produces denseness primarily around the global minimizers and not throughout the entire domain $\mathcal{D}$. However, achieving this in practice is very difficult, especially without prior knowledge of the global Lipschitz constant \cite{di2016direct, sergeyev2015deterministic}.
Hence, the stopping criterion for HALO is mainly based on the maximum number of function evaluations similarly as in \cite{liuzzi2010partition}, which is a commonly used criterion in practice.

However, we can show that the exact value of the global Lipschitz constant can be obtained when the algorithm iterates indefinitely to infinity.
First, we can highlight an important characteristic of the distance $\Delta_{i_k}$ in the next proposition
\begin{proposition}\label{prop:min_distance}
	Given the center $\mathbf{x}_{i_k}$ of its relative partition $\mathcal{D}_{i_k}$. If the algorithm satisfies Property \ref{property:nested}  and Property \ref{property:dist1}
	then the distance $\Delta_{i_k}$ from the nearest point of $\mathbf{x}_{i_k}$ satisfies the following limit
	\begin{equation}
		\lim\limits_{k\to+\infty}\Delta_{i_k} = 0
	\end{equation}
\begin{proof}
	If Property \ref{property:nested} holds the algorithm produces at least one sequence of strictly nested partitions $\{\mathcal{D}^{i_k}\}$ such that
	\begin{equation}\label{eq:limits_u_and_l}
    \lim_{k\to+\infty} ||\mathbf{u}_{i_k} - \mathbf{l}_{i_k}|| = 0
	\end{equation}
then we know that if the partition $\mathcal{D}_{i_k}$ is selected, then exactly $J=2 |\mathcal{P}_{i_{k}}|$ (see Eq. \ref{eq:set_P}) points are sampled inside $\mathcal{D}_{i_k}$. Consequently given Property \ref{property:dist1}, for every $j=1,\dots, J$ we have that
\begin{equation}
	\Delta_{i_k}=||\mathbf{x}_{i_k} - \mathbf{x}_{i_k}^j|| \leq ||\mathbf{u}_{i_k} - \mathbf{l}_{i_k}||
\end{equation}
 considering Eq. \ref{eq:limits_u_and_l} it follows that $\lim\limits_{k\to+\infty}\Delta_{i_k} = 0$
\end{proof}
\end{proposition}
Proposition \ref{prop:min_distance} states that the minimum distance $\Delta_{i_k}$ from the nearest point to $\mathbf{x}_{i_k}$ tends to zero as the iteration $k$ approaches infinity and as the partition $\mathcal{D}_{i_k}$ is continuously selected.
Thus, as a partition $\mathcal{D}_{i_k}$ is selected during the iterations, the information contained in the vector $\widetilde{\nabla}{f(\mathbf{x}_{i_k})}$ becomes more precise, so that we can introduce the following propositions
\begin{proposition}\label{prop:lim_gradient}
	Given the representative point $\mathbf{x}_{i_k}$ which is center of its relative partition $\mathcal{D}_{i_k}$. If HALO satisfies Property \ref{property:nested} and Property \ref{property:dist1} then the vector $\widetilde{\nabla}{f(\mathbf{x}_{i_k})}$ associated to the partition $\mathcal{D}_{i_k}$ converge to the true gradient at $\mathbf{x}_{i_k}$
	\begin{equation}
		\lim\limits_{k\to+\infty}\widetilde{\nabla}{f(\mathbf{x}_{i_k})} = \nabla{f(\mathbf{x}_{i_k})}
	\end{equation}
	\begin{proof}
		If Property \ref{property:dist1} is satisfied then there exist $k$ such that for every $\bar{k}>k$ the centroid $\mathbf{x}_{i_k}$ has its nearest point along all the coordinate directions. Then using the results from Property \ref{property:nested} and for the definition of derivative it follows that $\lim\limits_{k\to+\infty}\widetilde{\nabla}{f(\mathbf{x}_{i_k})} = \nabla{f(\mathbf{x}_{i_k})}$.
	\end{proof}
\end{proposition}
\begin{proposition}\label{prop:lim_L_global}
	If HALO satisfies Proposition \ref{prop:lim_gradient} then 
	\begin{equation}
		\lim_{k\to\infty} \max_{i_k \in \mathcal{I}_k}\{||\widetilde{\nabla}{f(\mathbf{x}_{i_k})}||\} = L
	\end{equation}
	\begin{proof}
		The proof is given directly by using the results of Proposition \ref{prop:lim_gradient} and applying it for every $i_k \in \mathcal{I}_k$ and taking the maximum.
	\end{proof}
\end{proposition}
Proposition \ref{prop:lim_gradient} further shows that as the selected partition becomes denser the associated vector $\widetilde{\nabla}{f(\mathbf{x}_{i_k})}$ becomes a more accurate approximation of the true gradient $\nabla{f(\mathbf{x}_{i_k})}$. As a result, from Proposition \ref{prop:lim_L_global} the global Lipschitz constant $L$ can be accurately estimated.

\subsection{A Detailed Implementation of HALO}\label{section:pseudocode_descr}
The simplified pseudocode of HALO is presented in Algorithm \ref{alg:halo}. In the following, we will describe the main steps of the pseudocode.
The core part of HALO resides within the \texttt{Main} function, which takes as input the termination criterion for the algorithm, either based on the maximum number of function evaluations or the maximum number of iterations. The user also provides the box constraints of the optimization problem $\mathcal{D}$.

Starting with the counter $k$ set to zero, lines 9 and 10 define the indices of the selected partitions $\mathcal{I}^{\star}_{k}$ and the matrix $\mathbf{S} \in \mathbb{R}^{|\mathcal{I}_{k}| \times N}$, where each entry contains the vector $\mathbf{s}$ defined in Eq. \ref{eq:sides}.
In line 11, a for loop begins, iterating over the indices of each selected partition. This loop is primarily responsible for computing $\Delta_{i^{*}_k}$ as described in Section \ref{subsubsection:DIRECT}. 
In line 17, another for loop starts, cycling through each element of the set $\mathcal{P}_{i^{*}_k}$ defined at line 14 and Eq. \ref{eq:set_P}. This loop handles the sampling process described in Section \ref{subsubsection:DIRECT}. Additionally, we compute the slopes by updating the matrix $\mathbf{G}\in \mathbb{R}^{|\mathcal{I}_{k}| \times N}$, which stores the gradient estimations discussed in Section \ref{subsubsection:adagrad}. The notation $\mathbf{G}_{i^{*}_k, p}$ indicates that the $p$th coordinate of the gradient relative to the $i^{*}_k$th partition is being updated.
After this inner for loop, line 28 performs a sorting of the indices of the set $\mathcal{T}_{i^{\star}_k}$ defined in Eq. \ref{eq:set_T}, and line 29 calls the function (line 68) responsible for dividing the set $\mathcal{D}_{i^{*}_k}$. 

Once the outer loop defined in line 11 is completed, based on the information in the matrix $\mathbf{G}$, we invoke the function that computes the values of the local Lipschitz constants, stored in the vector $\boldsymbol{\ell}\in \mathbb{R}^{|\mathcal{I}{k}|}$. Lines 31 and 32 define the vectors $\mathbf{v}\in \mathbb{R}^{|\mathcal{I}{k}|}$ and $\mathbf{h}\in \mathbb{R}^{|\mathcal{I}_{k}|}$, which contain half of the distances from the center of each partition $j$ to its vertices and the norm of each gradient estimation (i.e., the norm of each row of $\mathbf{G}$). Line 34 computes the estimation of the local Lipschitz constant in a vectorized form, following Eq. \ref{eq:formula_adaptive}.

The next task is to use the local Lipschitz constants collected in the vector $\boldsymbol{\ell}$ to compute the lower bounds and select the most promising partitions based on Def. \ref{def:selection}. The lower bound values are stored in the vector $\mathbf{r}\in \mathbb{R}^{|\mathcal{I}_{k}|}$ (line 37). The indices $q^{\star}_1$, $q^{\star}_2$, and $q^{\star}_3$ represent the index of the partition with the minimum lower bound (\ref{choice 2}), the partition with the minimum objective function value (\ref{choice 3}), and the partition with the largest diagonal that has the minimum lower bound (\ref{choice 1}). It is possible for one index to satisfy multiple conditions (e.g., $q^{\star}_1 = q^{\star}_2$), so care should be taken to avoid duplicate entries in the set $\mathcal{I}_k^{*}$.
In lines 42 and 47, we check if a local search has already been performed from the centroid of a partition and if it satisfies Eq. \ref{eq:local_beta}. If both conditions are met, the local search is initiated as explained in Section \ref{section:coupling strategy}. Note that the largest hyperrectangle $q^{\star}_3$ is always selected (line 52). Finally, the set of selected partitions is returned, and a new loop can begin, iterating over its elements.
\begin{figure}[htb!]
	\centering
	\begin{subfigure}[t]{1.\textwidth}
		\includegraphics[width=\textwidth]{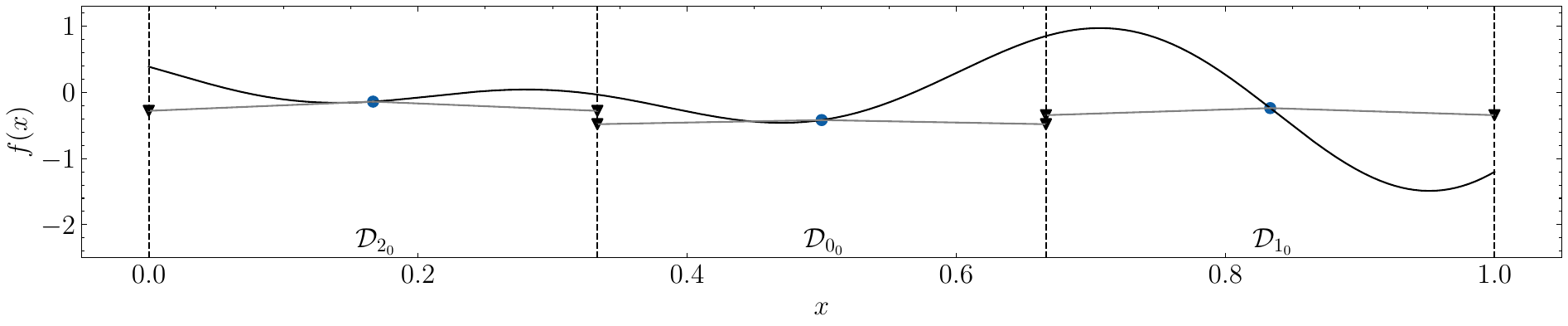}
		\caption{$k=0$}
		\label{fig:halo1d_k0}
	\end{subfigure}\\
	\begin{subfigure}[t]{1.\textwidth}
		\includegraphics[width=\textwidth]{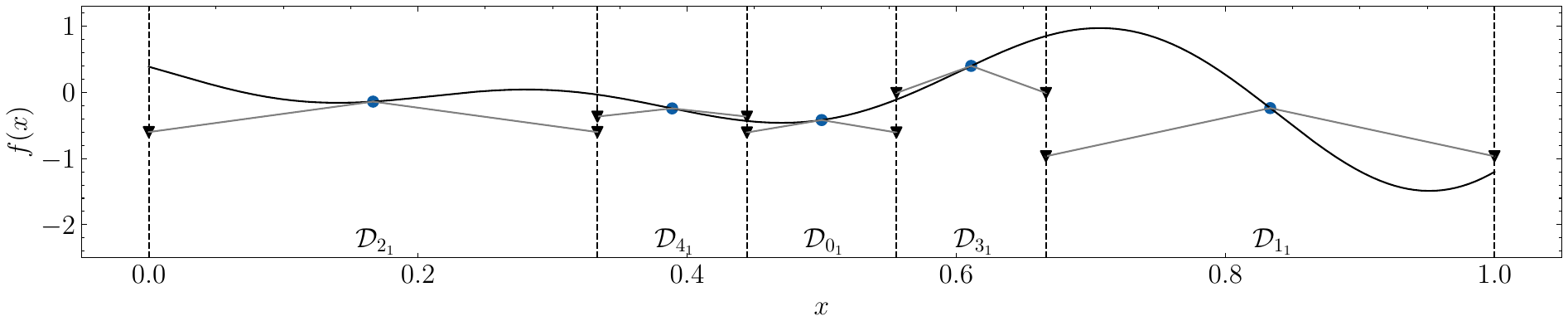}
		\caption{$k=1$}
	    \label{fig:halo1d_k1}
	\end{subfigure}\\
	\begin{subfigure}[t]{1.\textwidth}
		\includegraphics[width=\textwidth]{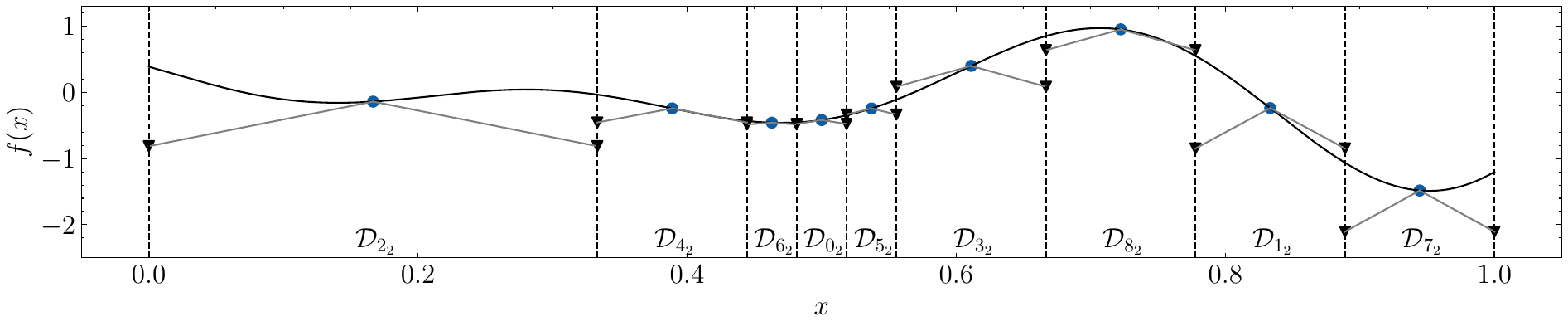}
		\caption{$k=2$}
		\label{fig:halo1d_k2}
	\end{subfigure}\\
	\begin{subfigure}[t]{1.\textwidth}
		\includegraphics[width=\textwidth]{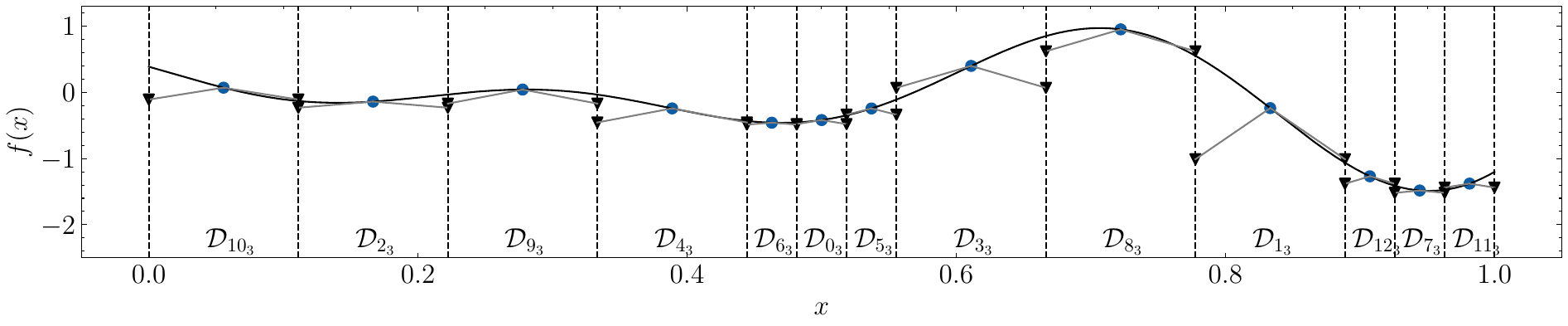}
		\caption{$k=3$}
		\label{fig:halo1d_k3}
	\end{subfigure}\\
	\caption{An example of HALO for a one-dimensional function. The blue rounded dots represent the centroids of each partition, the grey lines highlight the slope of the local Lipschitz constant, and the black triangles indicate the lower bound values at the boundaries of each partition.}
	\label{fig:halo1d}
\end{figure}
\subsubsection{Highlighting the Variable Importance and Problem Interpretability}\label{subsec:varimp}
As discussed in Section \ref{subsec:varimp}, HALO offers more than just an optimization routine; it also provides valuable insights into the importance of each variable in the black-box objective function. After HALO completes its search, we can compute the \emph{variable importance} using the matrix $\mathbf{G}$, which contains the approximate gradients $\widetilde{\nabla}{f(\mathbf{x}_{i_k})}$. By averaging these gradients as shown in Equation \ref{eq:VI}, we identify the directions in which the objective function is most sensitive.
In many practical applications, it is not only important to identify the global minimum but also to gain a deeper understanding of the variables that significantly impact the reduction of the objective function. Traditionally, such insights are obtained through sensitivity analysis, which involves computing gradients of the objective function. However, this process can be time-consuming, especially when dealing with computationally expensive simulations.

With HALO, we have an alternative approach to extracting similar information without the need for explicit and additional gradient computations. Once the optimization routine is completed, we can analyze the matrix $\mathbf{G}$ (introduced in Section \ref{section:pseudocode_descr}), which collects the vectors $\widetilde{\nabla}{f(\mathbf{x}_{i_k})}$ around each centroid $\mathbf{x}_{i_k}$. By examining this matrix, we can gain insights into the directions in which the objective function exhibits greater sensitivity during the iterations of HALO. We can simply determine the variable importance by averaging the rows $\widetilde{\nabla}{f(\mathbf{x}_{i_k})}$ of the matrix $\mathbf{G}$. This can be expressed as
\begin{equation}\label{eq:VI}
    \text{Variable Importance} = \frac{1}{|\mathcal{I}_k|}\sum_i^{|\mathcal{I}_k|} \widetilde{\nabla}{f(\mathbf{x}_{i_k})}
\end{equation}
we also apply a final normalization step so that the variable importance vector sums to one. 
\subsubsection{Example of HALO for Univariate Function}
To conclude this section, Fig. \ref{fig:halo1d} illustrates the behavior of HALO for a simple one-dimensional function. 
At iteration $k=0$, it can be observed that the values of the Lipschitz constants for the three partitions do not reflect the underlying trend of the objective function because these three points have similar objective function values, and there is a relatively large distance between them.
At iteration $k=1$, only the partition $\mathcal{D}_{0}$ is selected and divided because it satisfies \ref{choice 2}, \ref{choice 3}, and \ref{choice 1}. The lower bound values are already more precise at this stage.
At iteration $k=2$, two intervals are selected and divided: $\mathcal{D}_{0}$, which satisfies \ref{choice 3} (has the lowest function value), and $\mathcal{D}_{1}$, which satisfies \ref{choice 2} and \ref{choice 1} (it has the lowest overall lower bound and is the largest partition with the lowest lower bound).
Finally, at iteration $k=3$, two more intervals are divided: $\mathcal{D}_{2}$, which satisfies \ref{choice 1}, and $\mathcal{D}_{7}$, which satisfies \ref{choice 2} and \ref{choice 3}. The last figure (Fig. \ref{fig:halo1d_k3}), demonstrates that all the lower bound values are now quite precise. This entire process is automated in HALO, requiring no critical hyperparameters to be set by the user.

\section{Numerical Results}
In this section, we conduct a comprehensive evaluation of the HALO algorithm compared to other GO algorithms. The evaluation is performed on a diverse set of test functions organized into three different benchmarks. Additionally, we show how to extract insightful information from the objective function using our approach. 
Finally, we conduct a parametric study on the parameter $\beta$ from Eq. \ref{eq:local_beta} to provide general recommendations to users regarding its value.

The code for the algorithms and test functions used in this study can be found in the following GitHub repository \url{https://github.com/dannyzx/HALO}.

\subsection{Experimental Setup}

In the following paragraphs, we provide details about the algorithms, test functions, and stopping criteria used in our numerical experiments.

\subsubsection{Algorithms}

The following algorithms are included in our numerical experiments:
\begin{itemize}

\item $\textbf{HALO}_{\text{L-BFGS-B}}$ and $\textbf{HALO}_{\text{SDBOX}}$: These are two versions of HALO that are coupled with different local optimizers. $\text{HALO}_{\text{L-BFGS-B}}$ utilizes a quasi-Newton method for bound-constrained problems (L-BFGS-B) \cite{byrd1995limited}, while $\text{HALO}_{\text{SDBOX}}$ uses a derivative-free coordinate descent method with an Armijo-type linesearch for bound-constrained optimization \cite{lucidi2002-COA} (SD-BOX).
Both versions of HALO have a single hyperparameter, $\beta$, which is defined in Eq. \ref{eq:local_beta}. This parameter determines the minimum size a selected partition must have before a local search can be initiated from its centroid. For our experiments, we fix $\beta$ to $10^{-4}$. In $\text{HALO}_{\text{L-BFGS-B}}$, the gradients used in the Quasi-Newton code are computed using finite difference methods.

\item $\textbf{HLO}$: This is similar to $\text{HALO}_{\text{SDBOX}}$, but it does not use adaptive estimates of the local Lipschitz constants as defined in Eq. \ref{eq:formula_adaptive}. Instead, the lower bounds are calculated based on the estimate of the global Lipschitz constant only. In other words, all local Lipschitz constants are fixed as  $\tilde{L}_{i_k} = \max_{i_k \in \mathcal{I}k} ||\widetilde{\nabla}{f(\mathbf{x}_{i_k})}||$. This allows us to compare the convergence acceleration achieved by estimating the local Lipschitz constants using only the global Lipschitz constant estimate.
\item \textbf{DIRECT} \cite{jones1993-JOTA}: This is the classical DIRECT algorithm, with a fixed hyperparameter $\epsilon$ set to $10^{-4}|f_{\text{min}}|$, as suggested in \cite{jones1993-JOTA}.
\item \textbf{HPLOR} \cite{mockus2017application}: HPLOR is a Pareto-reduced Lipschitzian optimization method similar to the DIRECT algorithm, where only a subset of potentially optimal hyperrectangles is selected. Specifically, it chooses the one with the minimum function value and the one with the largest diagonal. This approach assumes the Lipschitz constant to be either very low or very high. To ensure a fair comparison, we introduced a hybridized version of the algorithm, referred to as Hybrid-PLOR (HPLOR), using the same coupling strategy as HALO and incorporating L-BFGS-B as the local optimizer.

\newpage
\scalebox{0.62}{\begin{minipage}{1.5\linewidth}
\begin{algorithm}[H]\label{algo:halo}
\scriptsize
  \DontPrintSemicolon
\Function{Main(\text{max\_iter}, \text{max\_fun\_eval}, $\mathcal{D}$)}{
     $k \assign 0$ \;
     $\mathcal{C}_k, \mathbf{G}, \mathbf{S}, \mathbf{f} \assign \emptyset, [[\,\,]], [[\, \,]], [\, \,]$ \;
    \For{$k$ \KwTo \text{max\_iter}}{
         \uIf{$k>0$}{
           $\boldsymbol{\ell} \assign$ \FuncCall{compute\_Lipschitz}{$\mathbf{G}$, $\mathbf{S}$}  \;
           $\mathcal{I}^{\star}_{k} \assign$ \FuncCall{selection}{$\mathcal{I}_{k}$, $\boldsymbol{\ell}$, $\mathcal{C}_k$} \;
        }
         \uElse{
           $\mathcal{I}^{\star}_{k} \assign \{ 0\} $\;
           $\mathbf{S}_{i_{k}} \assign \left[ |u_n - l_n|/2 \right]_{n=1}^{N} $\;
        }
    \For{$i^{\star}_k$ \KwIn $\mathcal{I}^{\star}_{k}$}{
        $\mathbf{s}_{i^{*}_{k}} \assign \mathbf{S}_{i^{*}_{k}}$\;
        $s_{i^{*}_{k}}^{\max} \assign \max\{\mathbf{s}_{i^{*}_{k}}\}$\;
        $\mathcal{P}_{i^{*}_k} \assign \argmax\{\mathbf{s}_{i^{*}_{k}}\}$\;
        $\Delta_{i^{*}_k} \assign \frac{2}{3} s_{i^{*}_{k}}^{\max}$\;
        $\mathcal{V}_{i^{\star}_k} \assign \emptyset$ \;
        \For{$p$ \KwIn $\mathcal{P}_{i^{*}_k}$}{
            $\mathbf{x}_{i^{*}_k}^{p_1} \assign \mathbf{x}_{i^{*}_k} + \Delta_{i^{*}_k} \mathbf{e}_{p}$ \;
            $f(\mathbf{x}_{i^{*}_k}^{p_1}) \assign $ \FuncCall{fun\_eval}{$\mathbf{x}_{i^{*}_k}^{p_1}$}\;
            $\mathbf{f}_{i^{*}_k + 1} \assign f(\mathbf{x}_{i^{*}_k}^{p_1})$\;
            $\mathbf{G}_{i^{*}_k + 1, p} \assign |f(\mathbf{x}_{i^{*}_k}^{p_1}) - f(\mathbf{x}_{i^{*}_k}))| / \Delta_{i^{*}_k} $\;
            $\mathbf{x}_{i^{*}_k}^{p_2} \assign \mathbf{x}_{i^{*}_k} - \Delta_{i^{*}_k} \mathbf{e}_{p}$ \;
            $f(\mathbf{x}_{i^{*}_k}^{p_2}) \assign $ \FuncCall{fun\_eval}{$\mathbf{x}_{i^{*}_k}^{p_2}$}\;
            $\mathbf{f}_{i^{*}_k + 2} \assign f(\mathbf{x}_{i^{*}_k}^{p_2})$\;
            $\mathbf{G}_{i^{*}_k + 2, p} \assign |f(\mathbf{x}_{i^{*}_k}^{p_2}) - f(\mathbf{x}_{i^{*}_k}))| / \Delta_{i^{*}_k} $\;
            $\mathbf{G}_{i^{*}_k, p} \assign |f(\mathbf{x}_{i^{*}_k}^{p_1}) - f(\mathbf{x}_{i^{*}_k}^{p_2}))| / 2\Delta_{i^{*}_k} $\;
            $\mathcal{T}_{i^{\star}_k} \assign \mathcal{T}_{i^{\star}_k} \cup \{\min \{f(\mathbf{x}_{i^{*}_k}^{p_1}), f(\mathbf{x}_{i^{*}_k}^{p_2})\}\}$
            }
        $\mathcal{U}_{i^{\star}_k} \assign \argsort \{\mathcal{T}_{i^{\star}_k}\}$ \;
        
        $\mathcal{I}_{k} \assign$ \FuncCall{partitioning}{$\mathcal{U}_{i^{\star}_k}, \mathbf{S}, \mathbf{V}$} \;
    }
    }
}
\Function{compute\_local\_Lipschitz($\mathbf{G}$)}{
    $\mathbf{v} \assign (||\mathbf{S}_j||)_{1\leq j \leq |\mathcal{I}_{k}|}$\;
    $\mathbf{h} \assign (||\mathbf{G}_j||)_{1\leq j \leq |\mathcal{I}_{k}|}$\;
    $\boldsymbol{\alpha} \assign \mathbf{v}/\sqrt{N}$\;
    $\boldsymbol{\ell} \assign \boldsymbol{\alpha} \odot\mathbf{h} + (1-  \boldsymbol{\alpha}) \odot \max\{\mathbf{h}\} $\;
    \Return $\boldsymbol{\ell}$\;
  }
\Function{selection($\mathcal{I}_{k}$, $\mathcal{C}_k$)}{
    $\mathbf{r} \assign \mathbf{f} - \mathbf{v} \odot \boldsymbol{\ell}$ \;
    $q^{\star}_1 \assign \argmin\{\mathbf{r}\}$ \;
    $q^{\star}_2 \assign \argmin\{\mathbf{f}\}$ \;
    $\mathcal{Q}_3 \assign \argmax\{\mathbf{v}\}$ \;
    $q^{\star}_3 = \argmin_{q_3 \in \mathcal{Q}_3} \{\boldsymbol{\ell}\}$ \;
    \uIf{$q^{\star}_1 \notin \mathcal{C}_k \, \KwAnd \, \mathbf{v}_{q^{\star}_1} \leq \beta$}{
        $\mathbf{x}^0 \assign \mathbf{x}_{q^{\star}_1}$ \;
        \FuncCall{local\_search}{$\mathbf{x}_{q^{\star}_1}, \mathcal{C}_k$}\;
    }
    \uElse{$\mathcal{I}^{\star}_{k} \assign \mathcal{I}^{\star}_{k} \cup \{q^{\star}_1\}$
    }
    \uIf{$q^{\star}_2 \notin \mathcal{C}_k \, \KwAnd \, \mathbf{v}_{q^{\star}_2} \leq \beta$}{
        $\mathbf{x}^0 \assign \mathbf{x}_{i^{*}_k}$ \;
        \FuncCall{local\_search}{$\mathbf{x}_{q^{\star}_2}, \mathcal{C}_k$}\;
    }
    \uElse{$\mathcal{I}^{\star}_{k} \assign \mathcal{I}^{\star}_{k} \cup \{q^{\star}_2\}$
    }
    $\mathcal{I}^{\star}_{k} \assign \mathcal{I}^{\star}_{k} \cup \{q^{\star}_3\}$\;
    \Return $\mathcal{I}^{\star}_{k}$\;
  }
\Function{local\_search($\mathbf{x}^0$, $\mathcal{C}_k$)}{
    $\mathbf{x}_{p_k} \assign \mathbf{x}^0$ \;
    \uIf{$\mathcal{C}_k = \emptyset$}{
    $\mathcal{C}_k \assign \{ j\in \mathcal{I}_k : ||\mathbf{x}_{p_k} - \mathbf{x}_{j_k}|| \leq r\}$ \;
    $\mathbf{x} \assign$ \FuncCall{local\_optimizer}{$\mathbf{x}_{p_k}$}\;
    \Return $\mathbf{x}$\;}
    \uElse{
    $\mathcal{B}_k \assign \{ j\in\mathcal{C}_k : ||\mathbf{x}_{p_k} - \mathbf{x}_{j_k}|| \leq r\}$ \;
    
    \uIf{$\mathcal{B}_k \neq \emptyset$}{
    $\mathcal{C}_k \assign \mathcal{C}_k \cup \{ j\in \mathcal{I}_k : ||\mathbf{x}_{p_k} - \mathbf{x}_{j_k}|| \leq r\}$  \;
    $\mathbf{x} \assign$ \FuncCall{local\_optimizer}{$\mathbf{x}_{p_k}$}\;
    \Return $\mathbf{x}$\;}

    \uElse{
        $\mathcal{C}_k \assign \mathcal{C}_k \cup p_k$ \;}
    }}
\Function{partitioning($\mathcal{U}_{i^{\star}_k},\mathcal{P}_{i^{\star}_k}, \mathbf{S}, \mathbf{V}$)}{
    \For{$u$ \KwIn $\mathcal{U}_{i^{\star}_k}$}{
            $m \assign \mathcal{P}_{i^{\star}_k, u} $ \;
            $\mathbf{S}_{i^{\star}_k, m} \assign \Delta/2$ \;
            $\mathbf{S}_{{i^{*}_k}+1} \assign \mathbf{S}_{i^{\star}_k}$ \;
            $\mathbf{S}_{{i^{*}_k}+2} \assign \mathbf{S}_{i^{\star}_k}$ \;
            $\mathcal{I}_{k} \assign \mathcal{I}_k \cup \{|\mathcal{I}_k| + 1, |\mathcal{I}_k| + 2\}$ \;
           }
    \Return $\mathcal{I}_{k}$\;

}
\caption{Hybrid Adaptive Lipschitzian Optimization (HALO)}
\label{alg:halo}
\end{algorithm}
\end{minipage}}

\item \textbf{DIRMIN} \cite{liuzzi2016exploiting}: an efficient hybridization of the DIRECT algorithm. It starts local optimizations from each potentially optimal hyperrectangle identified by DIRECT. We use the same derivative-free local optimizer as in $\text{HALO}_{\text{SDBOX}}$, which is the algorithm defined in \cite{lucidi2002-COA}. DIRMIN has shown improved convergence speed compared to DIRECT on various test functions \cite{liuzzi2016exploiting} and has shown to be among the best DIRECT-type approaches \cite{stripinis2022directgo, stripinis2023derivative, stripinis2024benchmarking}. It also outperformed other derivative-free hybrid GO algorithms in a real-world hull form shape optimization problem \cite{dagostino2018-AIAA-MAO}.

\item \textbf{CMA-ES} \cite{hansen2001completely}: a popular stochastic meta-heuristic approach known for its strong performance in various applications. It is considered one of the most competitive algorithms in its category \cite{hansen2010comparing}. We fix the hyperparameter $\sigma_0$ to be $1/4$ of the domain, as recommended by the author. Due to the stochastic nature of the algorithm, we run 20 independent experiments for each test function.

\item \textbf{L-SHADE} \cite{tanabe2014improving}: another popular stochastic method that improves upon the SHADE algorithm by incorporating linear population size reduction. L-SHADE has demonstrated strong performance in various benchmarks and optimization challenges.  We run 20 independent experiments to account for the stochastic nature of the algorithm and ensure reliable results. We used the implementation freely available at \footnote{\url{https://github.com/xKuZz/pyade}}.
\end{itemize}
\subsubsection{Test Functions}
We conducted our algorithm evaluations using three distinct benchmarks. Here are the details of each benchmark:
\begin{enumerate}
    \item \textbf{First Benchmark}: it is constructed using a function generator described in \cite{schoen1993wide}. We considered 100 test functions for each dimension $N$ ranging in $\{2,3, 4,6,8,10\}$, resulting in a total of 600 different test functions. The function generator requires two hyperparameters: the number of stationary points and their smoothness. For this experiment, we uniformly sampled the number of stationary points $S$ from the range $1$ to $100$, and the smoothness parameter was chosen uniformly from the range $[2, 3]$, following the approach in \cite{cassioli2013global}. 

\begin{figure}[htb!]
	\centering
	\begin{subfigure}[t]{0.3\textwidth}
		\includegraphics[width=\textwidth]{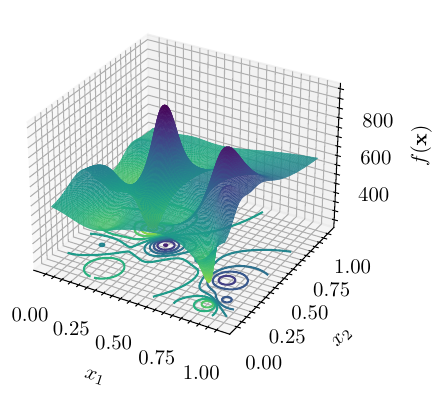}
		\caption{$S = 13$}
	\end{subfigure}
	\begin{subfigure}[t]{0.3\textwidth}
		\includegraphics[width=\textwidth]{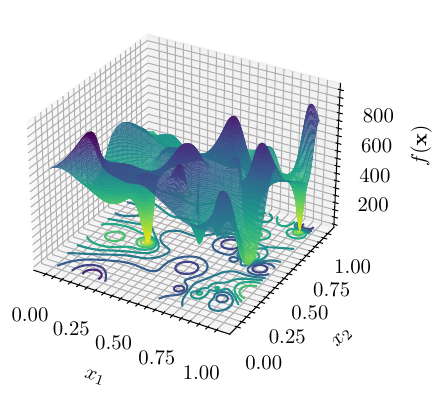}
		\caption{$S = 32$}
	\end{subfigure}
	\begin{subfigure}[t]{0.3\textwidth}
		\includegraphics[width=\textwidth]{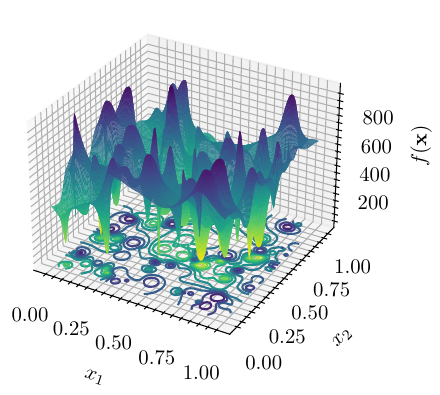}
		\caption{$S = 97$}
	\end{subfigure}
	\caption{Example of the first benchmark considered, $S$ represents the number of stationary points.}
\end{figure}
%
\item \textbf{Second Benchmark}: it consists of objective functions generated using the methodology described in \cite{gaviano2003algorithm}. The original implementation was in C programming language, but we rewrote it in Python for our experiments. Similarly to the first benchmark, we considered 100 test functions for each dimension $N$ ranging in $\{2,3, 4,6,8,10\}$, resulting in a total of 600 different test functions. This benchmark involves several hyperparameters: the distance $A$ from the global minimizer to the vertex of the paraboloid, the size of the basin of attraction $B$ of the global minimizer, and the number of local minima $C$. To create a challenging benchmark, we randomly set $A$ in the range $[0.8, 1)$, $B$ in the range $[0.1, 0.2)$, and $C$ from 3 to 10. These parameter settings ensure that the global minimizer is far from the vertex of the paraboloid and has a small basin of attraction.
%
\begin{figure}[htb!]
	\centering
	\begin{subfigure}[t]{0.3\textwidth}
		\includegraphics[width=\textwidth]{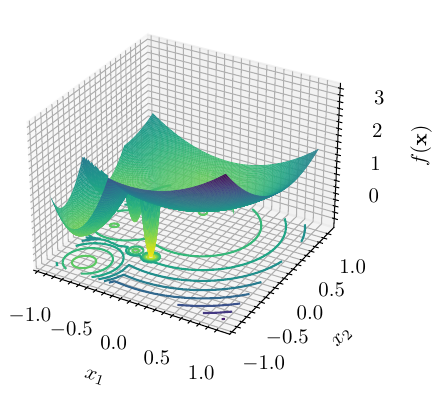}
		\caption{$A = 0.83, B=0.12,\\ C=7$}
	\end{subfigure}
	\begin{subfigure}[t]{0.3\textwidth}
		\includegraphics[width=\textwidth]{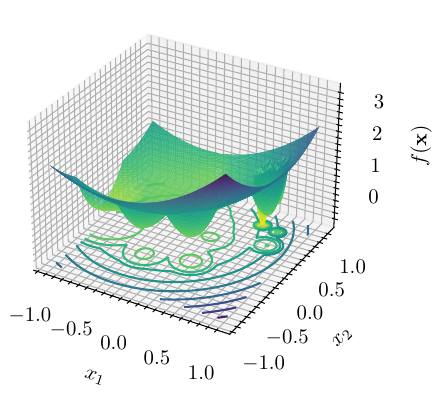}
		\caption{$A = 0.93, B=0.1,\\ C=8$}
	\end{subfigure}
	\begin{subfigure}[t]{0.3\textwidth}
		\includegraphics[width=\textwidth]{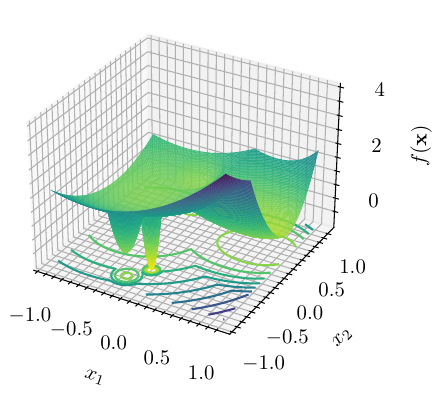}
		\caption{$A = 0.86, B=0.1,\\ C=4$}
	\end{subfigure}
	\caption{Example of the second benchmark considered, $A$ is the distance from the global minimizer to the vertex of the paraboloid, $B$ the radius of the basin of attraction of the global minimizer and $C$ the number of local minima.}
\end{figure}

    \item \textbf{Third Benchmark}: it comprises popular test functions frequently used to evaluate the performance of both local and global optimization algorithms, such as the Rosenbrock, Beale, Michalewicz, and Hartmann functions, as defined in \cite{jamil2013survey}. In total, we considered 191 different functions for this benchmark. When a test function is applicable to multiple dimensions $N$, we used the same function for different values of $N$ ranging in $\{2,3, 4,6,8,10\}$. Therefore, the total number of experiments using this benchmark is 496. In cases where a function has its global minimum at the center of the domain, we randomly shifted the global minimizer within the domain.
\begin{figure}[htb!]
	\centering
	\begin{subfigure}[t]{0.3\textwidth}
		\includegraphics[width=\textwidth]{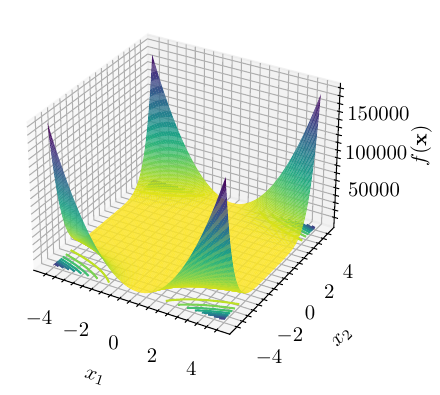}
		\caption{Beale.}
	\end{subfigure}
	\begin{subfigure}[t]{0.3\textwidth}
	\includegraphics[width=\textwidth]{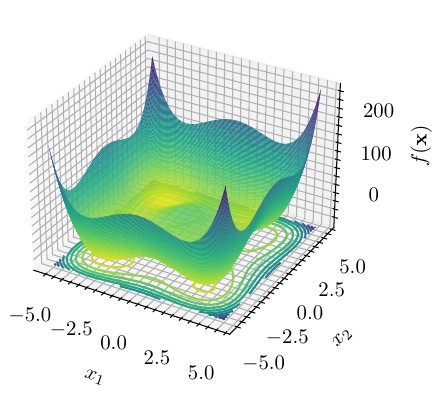}
	\caption{Styblinski Tang.}
	\end{subfigure}
	\begin{subfigure}[t]{0.3\textwidth}
	\includegraphics[width=\textwidth]{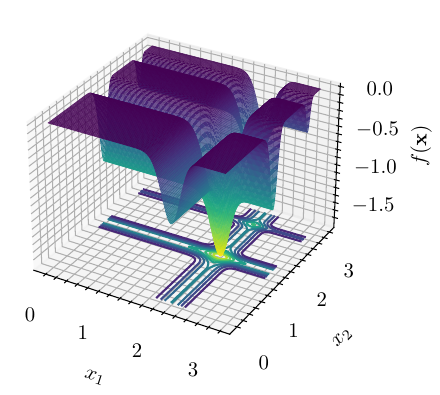}
	\caption{Michalewicz.}
	\end{subfigure}
	\caption{Example of the third benchmark considered.}
\end{figure}
\end{enumerate}
\subsubsection{Termination Criteria and Evaluation Metric}
The termination criteria for our experiments are defined based on two conditions. Firstly, the relative error to the global optimal objective function value, denoted as $f_{\text{glob}}$, is used as a stopping criterion. The termination condition is satisfied when the following inequality holds
\begin{equation}\label{eq:stop_crit}
\frac{f - f_{\text{glob}}}{|f_{\text{glob}}|} \leq 10^{-4}
\end{equation}
This condition ensures that the difference between the current objective function value ($f$) and the global optimal objective function value ($f_{\text{glob}}$) is within a small relative error tolerance of $10^{-4}$.

The second termination criterion is based on the maximum number of function evaluations. In our experiments, we fixed this limit to 50000, irrespective of the problem's dimensionality. If an algorithm exceeds 50000 function evaluations without satisfying the criterion defined in Eq. \ref{eq:stop_crit}, the run is considered to have failed.

These termination criteria provide a balance between achieving accurate solutions and limiting the computational effort. In particular, the budget of 50,000 evaluations is generally considered sufficient in the context of real-world black-box problems, as using significantly more would be unrealistic.

To evaluate the performance of the algorithms given the presence of large number of objective functions in each benchmark, we will provide a summary of the results using operational characteristics \cite{strongin2013global}, also known as performance profiles \cite{dolan2002benchmarking}. 
More in particular the operational characteristic is defined as follows
\begin{equation}
    c(\gamma) = \frac{\text{Number of problems solved in less than $\gamma$ function evaluations}}{\text{Total number of problems}} 
\end{equation}
This approach allows us to effectively summarize and compare the performance of different algorithms across multiple test functions. To provide also a metric to evaluate algorithms we can compute the area under the operational characteristic (AUOC) 
\begin{equation}
    \text{AUOC} = \frac{1}{\gamma_{\text{max}}} \int_{0}^{\gamma_{\text{max}}} c(\gamma) \, d\gamma
\end{equation}
which is properly normalized so that AUOC $\in [0, 1]$. 
We will also present the percentage of problems solved, and an average number of function evaluations to further assess algorithm performance.
\subsection{Comparison Results and Discussion}
The numerical results highlight the significant advantage of using estimates of the local Lipschitz constant for accelerating convergence to the global optimum. Algorithms like $\text{HALO}_{\text{L-BFGS-B}}$ and $\text{HALO}_{\text{SDBOX}}$, which incorporate local Lipschitz estimates, outperformed direct-based approaches and those relying on sampling potentially optimal hyperrectangles. These methods consistently achieved higher total AUOC values, demonstrating their ability to more efficiently solve complex optimization problems.

In comparison, popular metaheuristics such as CMA-ES and L-SHADE showed weaker performance, failing to match the effectiveness of the Lipschitz-based algorithms. Despite their established use in global optimization, these methods did not achieve the same level of consistency or efficiency across the benchmarks.

Overall, the results confirm that incorporating local Lipschitz constant estimates can significantly improve optimization performance, offering a clear advantage over traditional DIRECT-type and metaheuristic approaches.
\begin{figure}[htb!]
	\centering
	
		\includegraphics[width=1.0\textwidth]{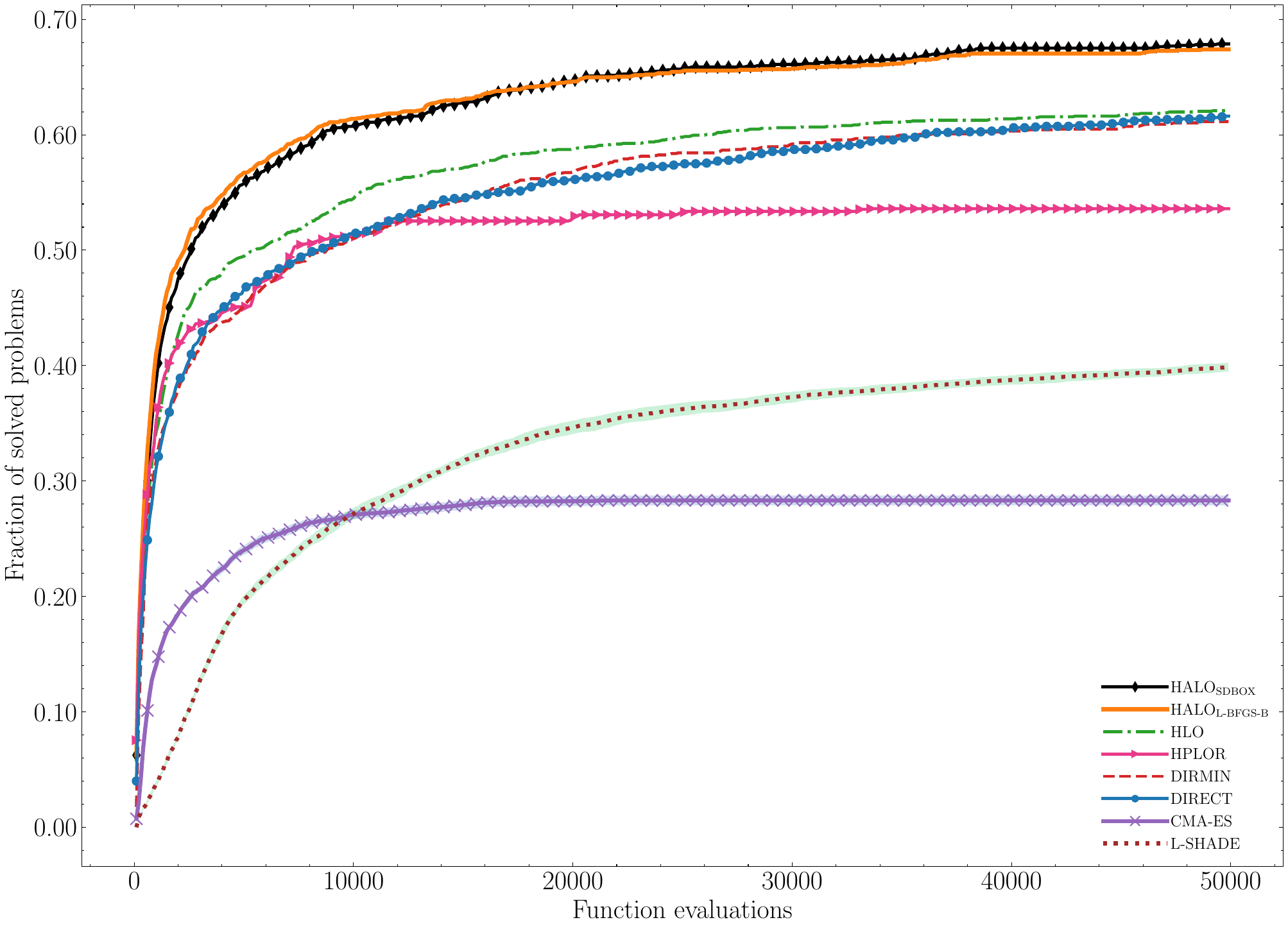}
		\caption{Operational characteristics for all the test functions considered in this study.}
		\label{fig:oc_all}
\end{figure}
\begin{table}[htb!]
\caption{Summary of the numerical results for different optimization algorithms across three benchmark tests. The table presents the total AUOC for each algorithm, along with their performance in the first, second, and third benchmarks.}
	\begin{tabular}{lcccc}\toprule
		\multicolumn{5}{c}{Numerical results} \\
		\midrule
		Algorithm & Total AUOC & \multicolumn{3}{c}{AUOC} \\ 
		\cmidrule(lr){3-5}
		         &  & First benchmark & Second benchmark & Third benchmark \\ 
		\midrule
	   $\text{HALO}_{\text{L-BFGS-B}}$ & $0.631$ & $0.773$ & $0.425$ & $0.709$ \\
	    $\text{HALO}_{\text{SDBOX}}$ & $0.630$ & $0.768$ & $0.427$ & $0.710$ \\
	    HLO     & $0.573$ & $0.620$ & $0.412$ & $0.710$ \\
        HPLOR	& $0.515$ & $0.520$ & $0.356$ & $0.701$ \\  
		$\text{DIRMIN}$	& $0.552$ & $0.759$ & $0.218$ & $0.707$ \\  
		DIRECT & $0.552$ & $0.684$ & $0.356$ & $0.629$ \\  
		CMA-ES & $0.270$ & $0.220$ & $0.0$ & $0.655$ \\  
  		L-SHADE & $0.324$ & $0.317$ & $0.103$ & $0.600$ \\  
		\bottomrule
	\end{tabular}
	\label{table:res_all} 
\end{table}
\subsubsection{First Benchmark Results}
According to Table \ref{table:res_all} the best-performing algorithm is $\text{HALO}_{\text{L-BFGS-B}}$, while the CMA-ES algorithm appears to be the least effective.

Notably, there is no significant difference in performance between the two versions of HALO equipped with different local optimizers. However, an interesting observation is a difference in performance between $\text{HALO}_{\text{SDBOX}}$ and HLO. $\text{HALO}_{\text{SDBOX}}$ solves more problems than HLO and demonstrates a much faster convergence, as evident from the operational characteristic plot shown in Fig. \ref{fig:oc_2}. Additionally, HLO performs worse than DIRECT. Upon analyzing this discrepancy, we found that HLO often selects only two hyperrectangles based on the estimate of the global Lipschitz constant. These hyperrectangles have the minimum objective function value (selected apriori) and the hyperrectangle with the largest diagonal and lowest objective function value. Consequently, if the global minimum is not within these two partitions, the overall convergence of the HLO algorithm is compromised. The same is valid for the HPLOR algorithm.

\begin{figure}[htb!]
	\begin{subfigure}[t]{0.32\textwidth}
    	\includegraphics[width=1\textwidth]{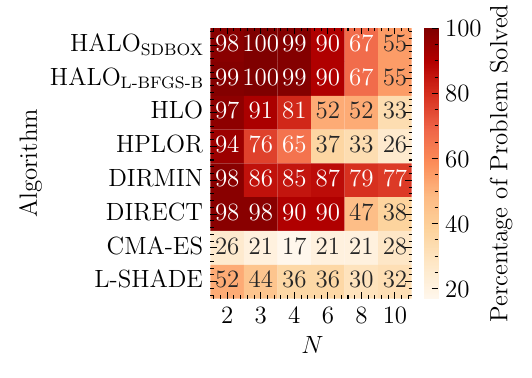}
    	\caption{First benchmark.}
    	\label{fig:dim_2}
	\end{subfigure}
	\begin{subfigure}[t]{0.32\textwidth}
    	\includegraphics[width=1\textwidth]{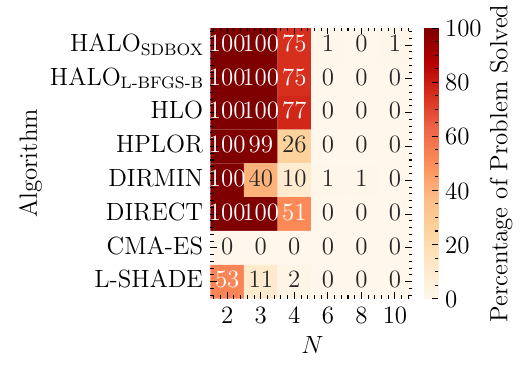}
    	\caption{Second benchmark.}
    	\label{fig:dim_3}
	\end{subfigure}
	\begin{subfigure}[t]{0.32\textwidth}
		\includegraphics[width=1\textwidth]{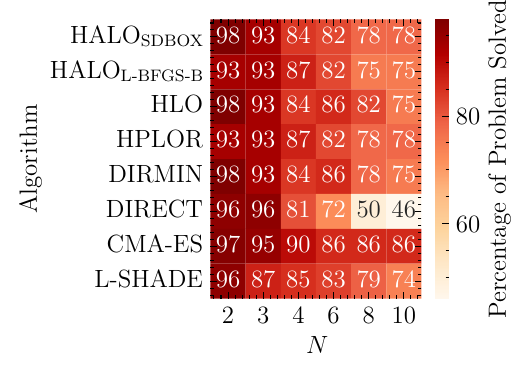}
		\caption{Third benchmark.}
		\label{fig:dim_1}
	\end{subfigure}
	\caption{Percentage of the problem solved varying the dimensionality $N$ for the three benchmarks.}
\end{figure}

In this benchmark, DIRMIN significantly accelerates the convergence of DIRECT. Particularly interesting is the performance of DIRMIN as the problem's dimensionality varies, as shown in Fig. \ref{fig:dim_2}. It can be observed that for $N \geq 6$, DIRMIN outperforms DIRECT, and for $N \geq 8$, it outperforms all the algorithms considered in this benchmark. 
This highlights that DIRMIN is less affected by the curse of dimensionality, possibly due to its effective utilization of the local optimization algorithm, which enables it to capture the local trends of the objective function.

Lastly, it is worth noting that CMA-ES often gets trapped in bad local minima for the majority of the objective functions, regardless of the problem's dimensionality (see Fig. \ref{fig:dim_2}). Similar behavior for L-SHADE although it performs better than CMA-ES.
%
\subsubsection{Second Benchmark Results}
In this benchmark, the best performing algorithm is $\text{HALO}_{\text{SD-BOX}}$, as it achieves the highest AUOC and solves the highest percentage of problems (see Table \ref{table:res_all}. However, its performance is quite similar to $\text{HALO}_{\text{SD-BOX}}$. In contrast, the CMA-ES algorithm remains the least effective in this benchmark.

Our analysis indicates that this highly challenging benchmark evaluates the capability of a global optimization algorithm to explore the entire search space efficiently. An algorithm that overly emphasizes exploitation rather than exploration could struggle with many of the test functions in this benchmark. This is due to the small basin of attraction for the global minimum and the fixed relatively high distance from the vertex of the paraboloid (acting as a 'trap') to the global minimizer in all the test functions considered.

It is also interesting to analyze the difference between using the local Lipschitz constant, as in $\text{HALO}_{\text{SDBOX}}$, or the global Lipschitz constant, like in HLO. Additionally, this benchmark confirms that using an estimate of the local Lipschitz constant for each partition leads to better performance. $\text{HALO}_{\text{SDBOX}}$ solves more problems than HLO and uses fewer function evaluations.

We can also observe that HPLOR, which does not use any estimate of the Lipschitz constants, performs worse compared to HALO and HLO. This suggests that leveraging information about the Lipschitz constants, both local (HALO) and global (HLO), provides valuable insights that help accelerate convergence toward the global optimum for this set of functions

When comparing the performance of DIRECT and its hybrid counterpart DIRMIN, it is evident that DIRMIN is significantly less effective than DIRECT. This is primarily because too many local searches start and end up at the vertex of the paraboloid.

It appears that CMA-ES encounters significant difficulty in solving any of the objective functions. Our results indicate that almost $100 \%$ of the time, CMA-ES gets trapped at the vertex of the paraboloid. 
Also L-SHADE experiences some difficulty with this benchmark, as is only able to solve $11.4$ \% of the problems.

It should be noted that the GKLS generator was reimplemented in this work, and since the full parameter sets used in prior studies were not released, the generated test functions cannot be replicated exactly, making direct comparison with previous results infeasible.

\subsubsection{Third Benchmark Results}
The numerical results for the third benchmark indicate that $\text{HALO}_{\text{SD-BOX}}$ and $\text{HLO}$ are the best-performing algorithm, according to Table \ref{table:res_all}. Conversely, the L-SHADE algorithm appears to be relatively less effective.

\begin{figure}[htb!]
	\centering
	\begin{subfigure}[t]{0.32\textwidth}
		\includegraphics[width=1.0\textwidth]{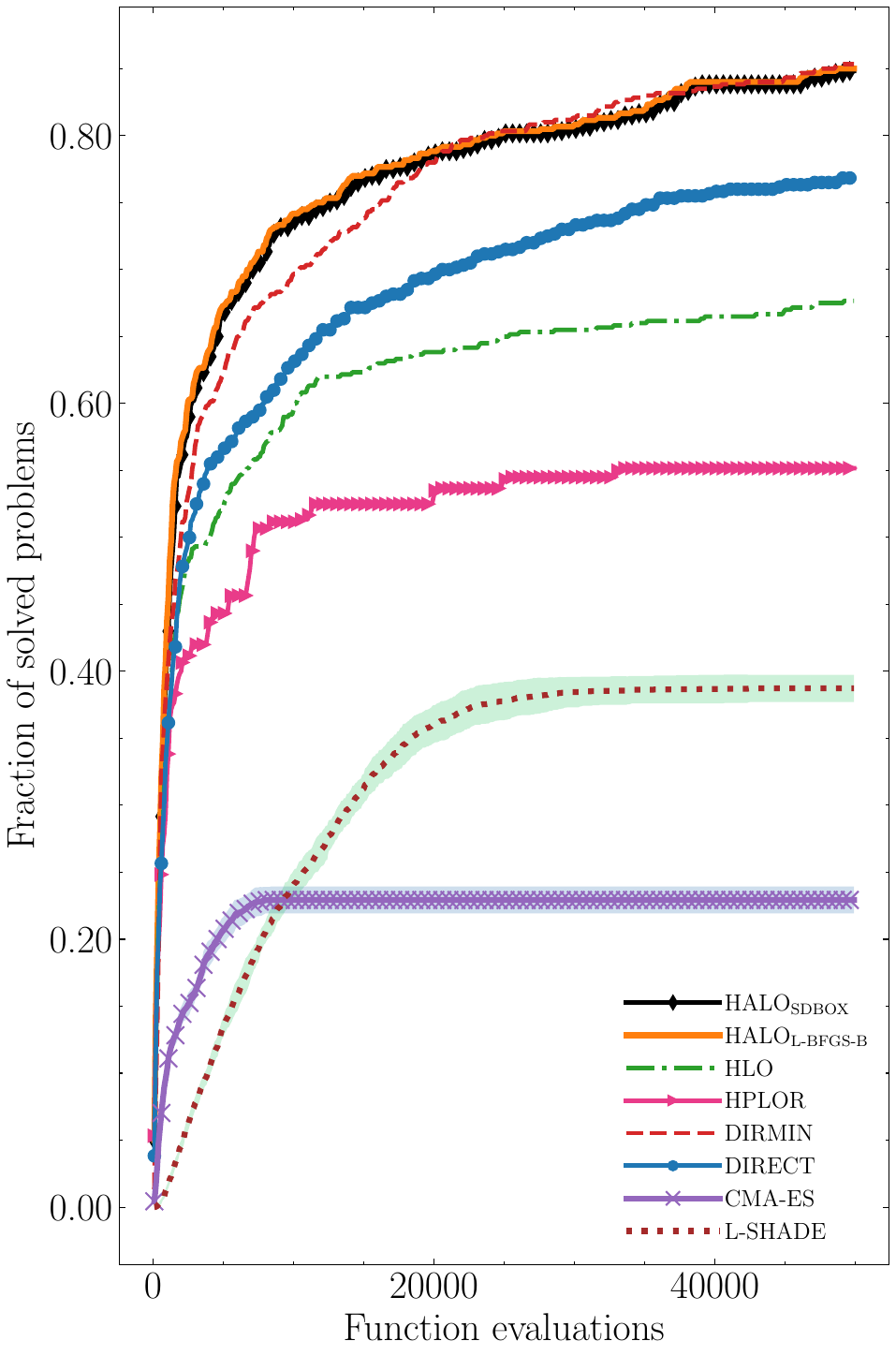}
		\caption{First benchmark.}
		\label{fig:oc_2}
	\end{subfigure}
	\begin{subfigure}[t]{0.32\textwidth}
		\includegraphics[width=1.0\textwidth]{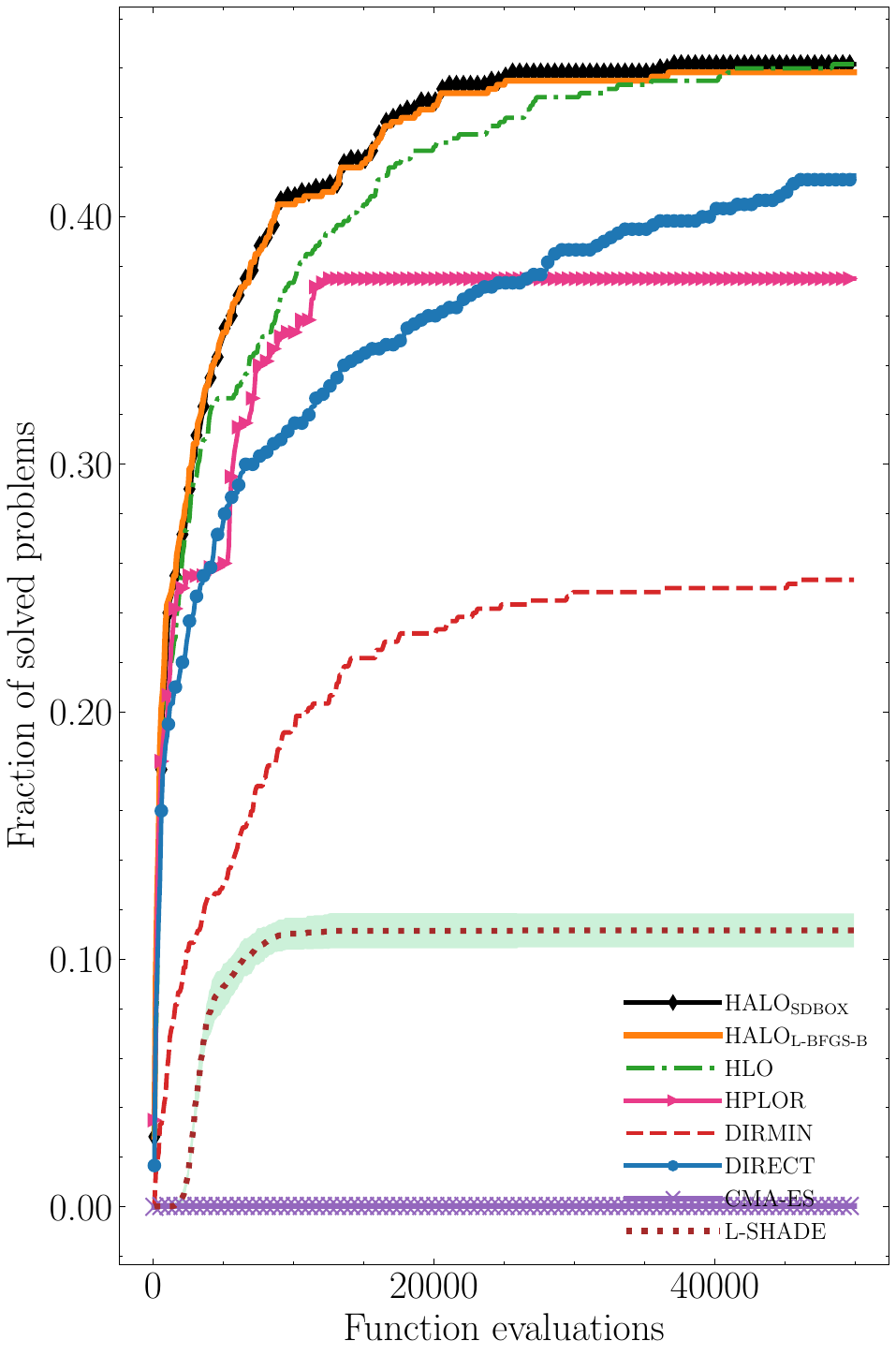}
		\caption{Second Benchmark.}
		\label{fig:oc_3}
	\end{subfigure}
	\begin{subfigure}[t]{0.32\textwidth}
		\includegraphics[width=1.0\textwidth]{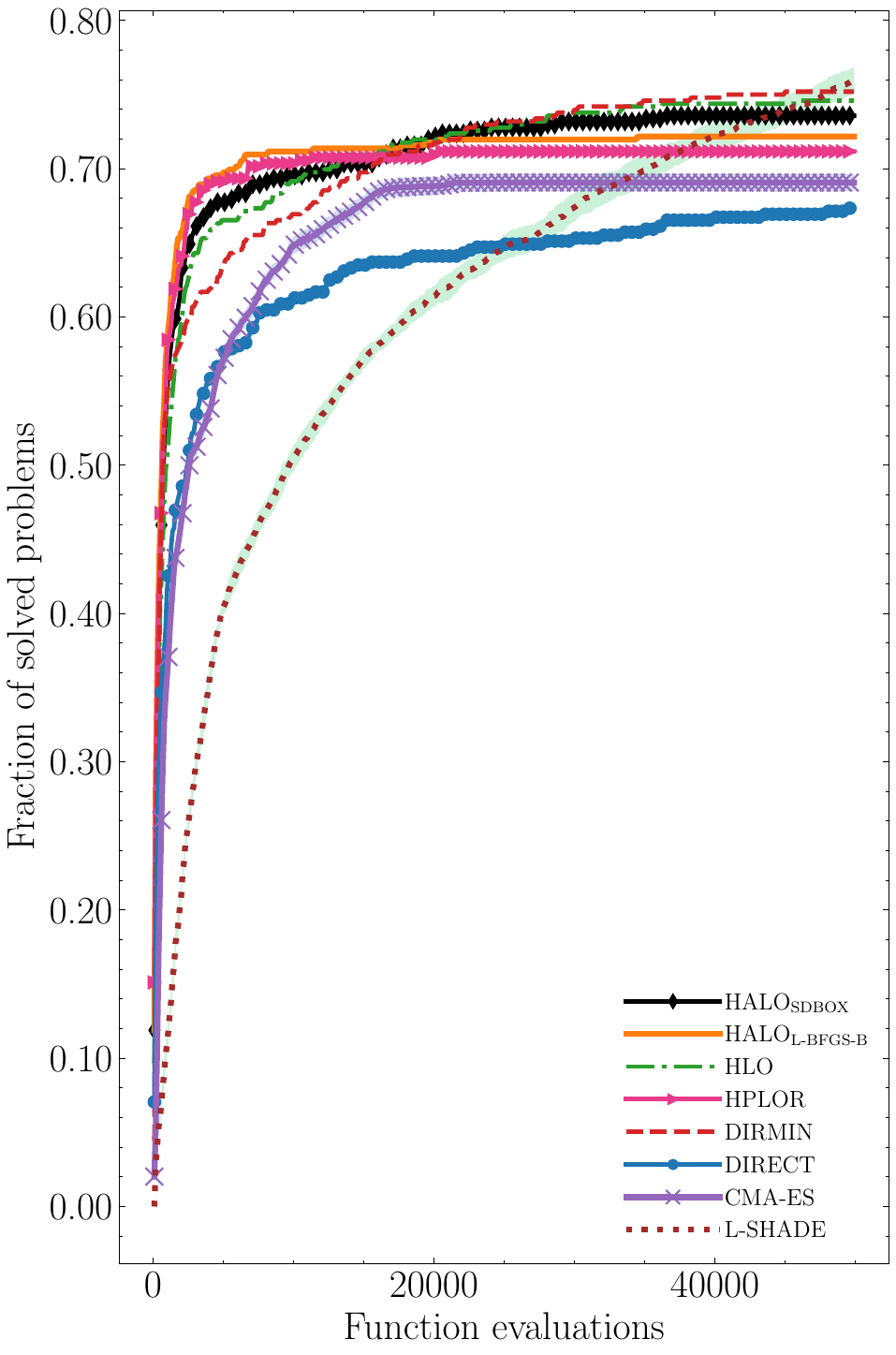}
		\caption{Third benchmark.}
		\label{fig:oc_1}
	\end{subfigure}
	\caption{Operational characteristics for the three benchmarks.}
\end{figure}

\begin{figure}[htb!]
	\begin{subfigure}[t]{0.32\textwidth}
    	\includegraphics[width=1\textwidth]{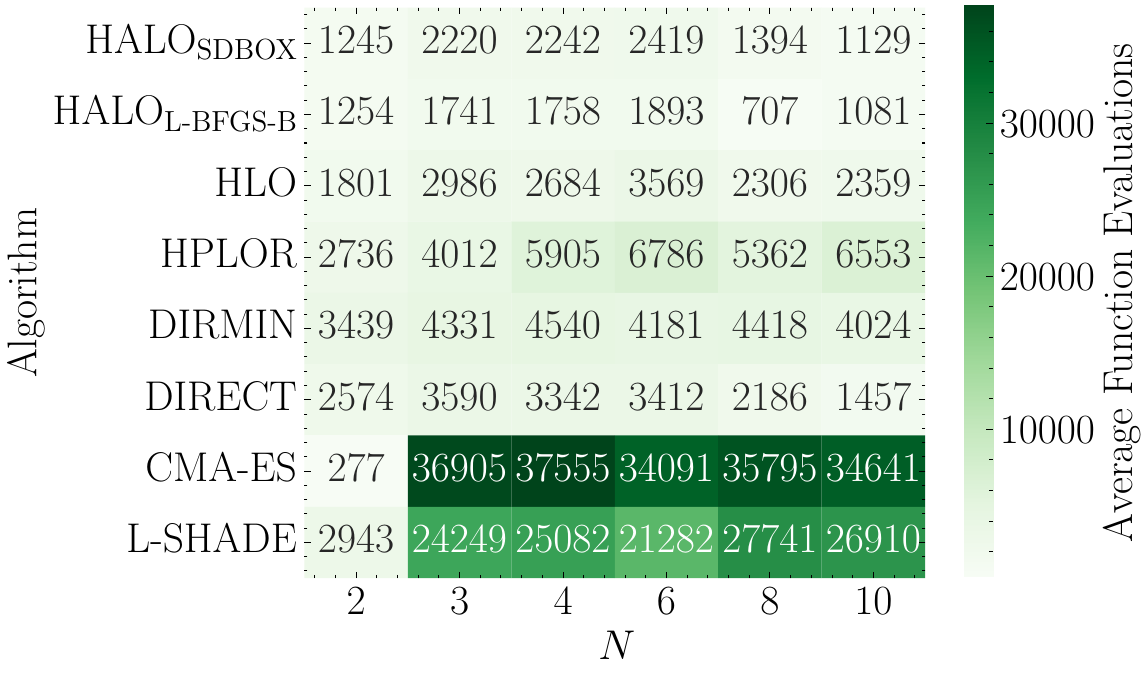}
    	\caption{First benchmark.}
    	\label{fig:dim_fe_2}
	\end{subfigure}
	\begin{subfigure}[t]{0.32\textwidth}
    	\includegraphics[width=1\textwidth]{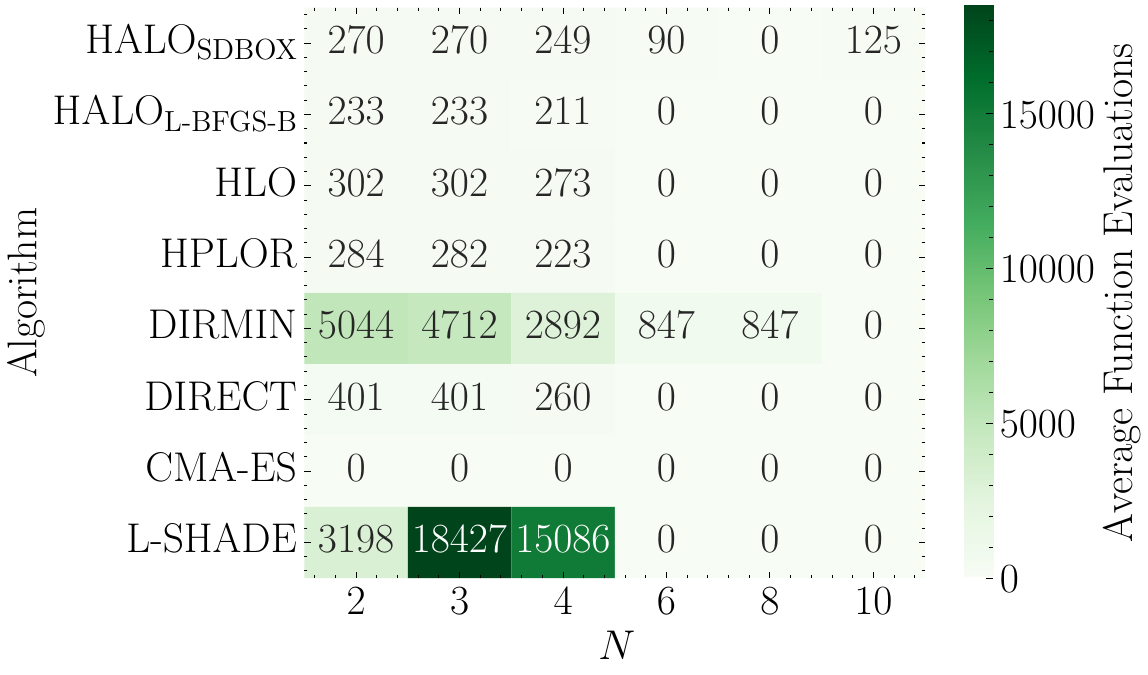}
    	\caption{Second benchmark.}
    	\label{fig:dim_fe_3}
	\end{subfigure}
	\begin{subfigure}[t]{0.32\textwidth}
		\includegraphics[width=1\textwidth]{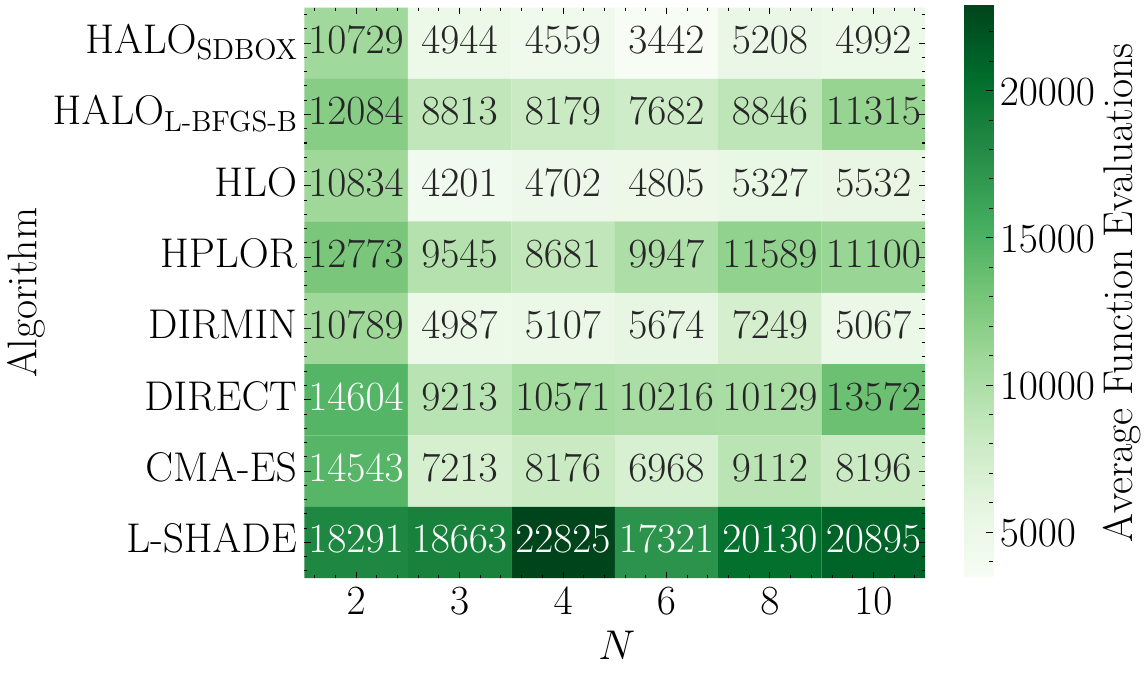}
		\caption{Third benchmark.}
		\label{fig:dim_fe_1}
	\end{subfigure}
	\caption{Average number of functions evaluations varying the dimensionality $N$ for the three benchmarks.}
\end{figure}

While $\text{HALO}_{\text{SDBOX}}$ and HLO exhibit similar problem-solving capabilities, $\text{HALO}_{\text{SDBOX}}$ achieves the solutions with fewer function evaluations than HLO (see Fig. \ref{fig:dim_fe_1}. This confirms that, in this benchmark as well, utilizing an estimate of the local Lipschitz constant accelerates the algorithm’s convergence towards the global minimizers more effectively than relying solely on the global Lipschitz constant. This observation is supported by the operational characteristics shown in Fig. \ref{fig:oc_1}.
Although DIRMIN significantly accelerates the convergence of DIRECT, it does not outperform HALO and HLO in terms of AUOC and function evaluations. This is likely because DIRECT-type approaches, based on potentially optimal hyperrectangles, tend to sample many of them, leading to extensive exploration of sub-optimal regions in the domain $\mathcal{D}$.
CMA-ES performs better than the DIRECT algorithm but remains less competitive compared to the other approaches considered. 
The behavior of L-SHADE is particularly notable when observing the operational characteristic in Fig. \ref{fig:oc_1}. The algorithm requires a large number of function evaluations to solve a substantial percentage of problems, which explains the low AUOC value it obtained.

Fig. \ref{fig:dim_1} highlights the performance of all algorithms across varying dimensionalities $N$. It can be observed that $\text{HALO}_{\text{L-BFGS-B}}$, $\text{HALO}_{\text{SDBOX}}$, HLO, and DIRMIN exhibit similar performance across different dimensionalities. CMA-ES also shows similar performance, slightly underperforming for $N = 2$ but remaining competitive for other values of $N$. On the other hand, DIRECT demonstrates a marked performance deterioration for $N \geq 6$ in this benchmark.

To further analyze the algorithm performances, we divided the functions in this benchmark into two groups. The first group consists of functions exhibiting a simple trend that attracts the algorithm towards a global minimizer. The second group contains functions with more complex behavior, making it more challenging to locate the global minimizer. The 'simple' group consists of 77 functions, while the 'hard' group consists of 116 functions.

\begin{figure}[htb!]
	\centering
	\begin{subfigure}[t]{0.48\textwidth}
		\includegraphics[width=1.0\textwidth]{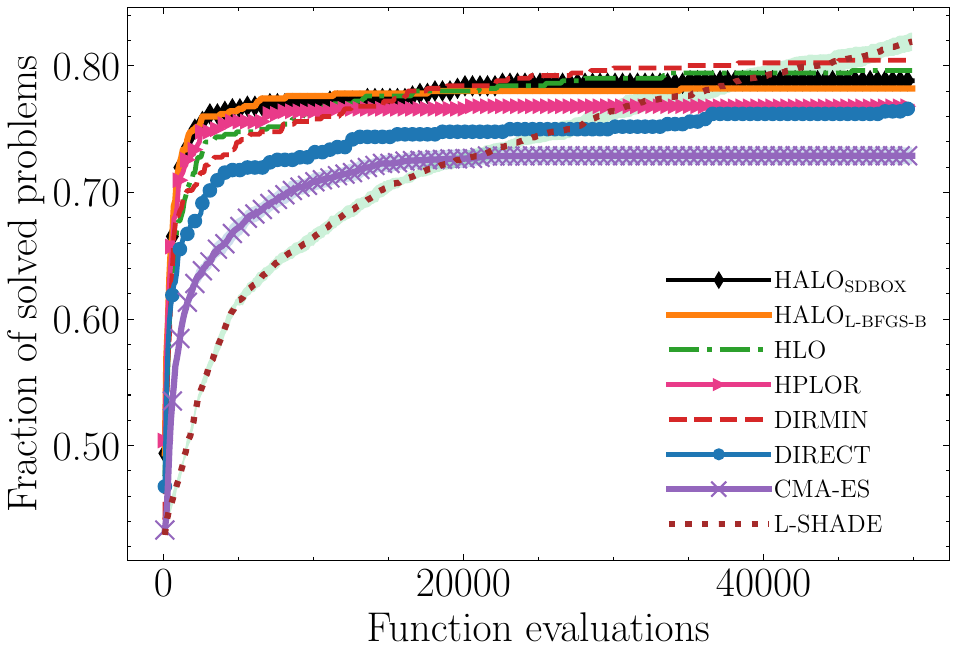}
		\caption{Hard functions.}
		\label{fig:scipy_hard_oc}
	\end{subfigure}
	\begin{subfigure}[t]{0.48\textwidth}
		\includegraphics[width=1.0\textwidth]{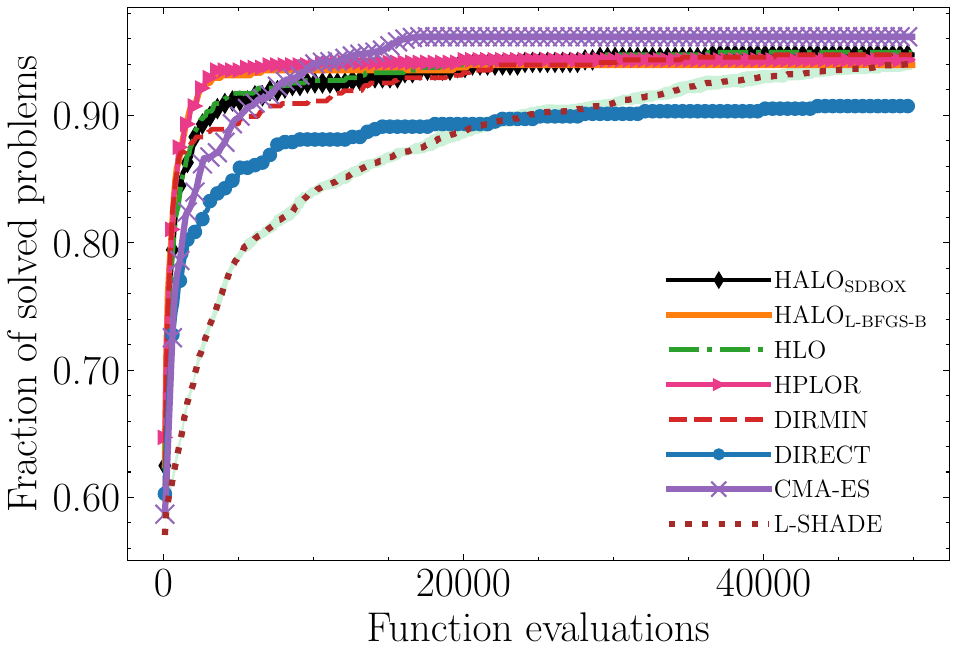}
		\caption{Simple functions.}
		\label{fig:scipy_simple_oc}
	\end{subfigure}
	\caption{Operational characteristics for the two groups (hard and simple) of functions belonging to the third benchmark.}
\end{figure}

We can make several observations by analyzing the operational characteristics in Fig. \ref{fig:scipy_hard_oc} and Fig. \ref{fig:scipy_simple_oc}. Firstly, CMA-ES performs poorly on the hard functions but performs remarkably on the simple functions group. This suggests that CMA-ES focuses more on exploitation than exploration. Another interesting finding is the comparison between HLO and $\text{HALO}_{\text{SDBOX}}$. In the simple group, their performances are similar, but $\text{HALO}_{\text{SDBOX}}$ outperforms HLO in the hard functions group. This indicates that when the objective function exhibits simple behavior, the information provided by the local Lipschitz constants used in $\text{HALO}_{\text{SDBOX}}$ does not significantly improve performance compared to using only the global Lipschitz constant in HLO. However, when the objective function displays more chaotic behavior, the information from the local Lipschitz constant becomes highly valuable. Regarding the comparison between $\text{HALO}_{\text{L-BFGS-B}}$ and $\text{HALO}_{\text{SDBOX}}$, as seen in Fig. \ref{fig:scipy_hard_oc} and Fig. \ref{fig:scipy_simple_oc}, the performance difference is minimal in the hard function group, while $\text{HALO}_{\text{L-BFGS-B}}$ outperforms $\text{HALO}_{\text{SDBOX}}$ in the simple function group. This suggests that the gradient-based local optimizer (L-BFGS-B) is more efficient than the derivative-free optimizer (SDBOX) when the objective function exhibits simple behavior. 

\subsection{Assessing Variable Importance and Interpretability with HALO}\label{subsec:var_selec}

As discussed in Section \ref{subsec:varimp}, HALO offers more than just an optimization routine; it also provides valuable insights into the importance of each variable in the objective function. After HALO completes its search, we can compute the \emph{variable importance} using the matrix $\mathbf{G}$, which contains the approximate gradients $\widetilde{\nabla}{f(\mathbf{x}_{i_k})}$. By averaging these gradients as shown in Equation \ref{eq:VI}, we identify the directions in which the objective function is most sensitive.

\begin{figure}[htb!]
	\centering
	\begin{subfigure}[t]{0.32\textwidth}
	\includegraphics[width=\textwidth]{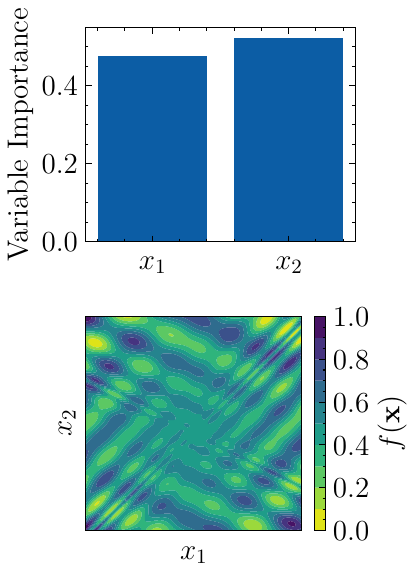}
	\caption{Egg Holder.}
        \label{fig:eggholder_2d}
	\end{subfigure}
	\begin{subfigure}[t]{0.325\textwidth}
	\includegraphics[width=\textwidth]{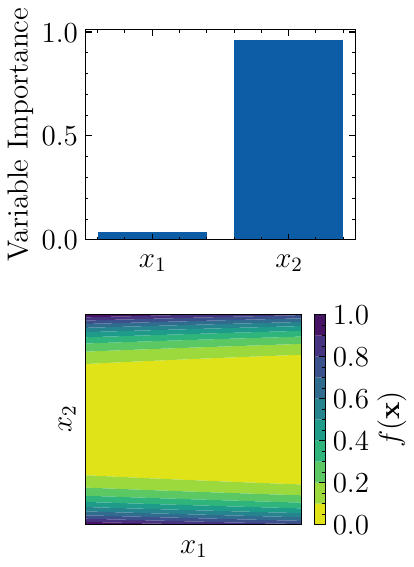}
	\caption{Dixon Price.}
     \label{fig:dixonprice_2d}
	\end{subfigure}
	\begin{subfigure}[t]{0.32\textwidth} 
	\includegraphics[width=\textwidth]{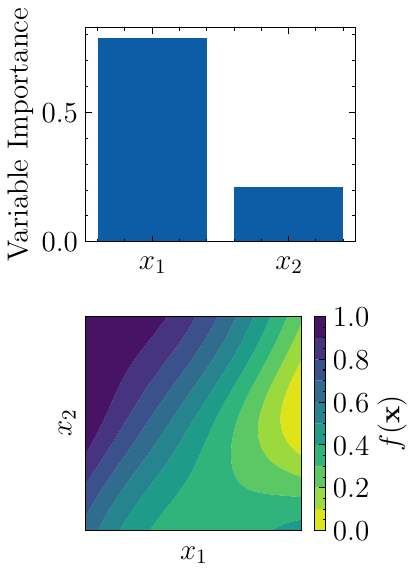}
	\caption{Adjiman.}
        \label{fig:Adjiman_2d}
	\end{subfigure} 
 \caption{Summary of the variable importance computed for various objective functions, extracted from HALO at the end of its iterations. The first row of figures shows the variable importance while in the second row the related objective function.}
 \label{fig:VI}
\end{figure}

In Fig. \ref{fig:VI}, we provide an illustrative example of the variable importance for three different objective functions: Eggholder, Dixon Price, and Adjiman respectively in Fig. \ref{fig:eggholder_2d}, \ref{fig:dixonprice_2d}, and \ref{fig:Adjiman_2d}.
By examining the variable importance defined in Eq. \ref{eq:VI}, we can observe that it accurately identifies the directions in which the objective function exhibits the most significant variation.

This information can be valuable in understanding the underlying dynamics of the objective function and identifying the variables that have the most significant influence on its behavior. 
This knowledge not only enhances our understanding of the problem but can also guide further investigations, inform decision-making processes, and potentially lead to improvements in the optimization problem formulation or solution strategies.
%
%

Compared to traditional sensitivity analysis techniques, leveraging the information available in $\mathbf{G}$ from HALO offers several advantages. First, it eliminates the need for additional gradient computations, saving computational time and resources, particularly when dealing with computationally expensive simulations. Second, it provides insights specific to the optimization process, capturing the dynamics and patterns observed during the search for the global minimum. 
\subsection{Hyperparameter Sensitivity Analysis}\label{subsec:PSA}
\begin{figure}[htb!]
	\centering
	\begin{subfigure}[t]{0.48\textwidth}
		\includegraphics[width=1.0\textwidth]{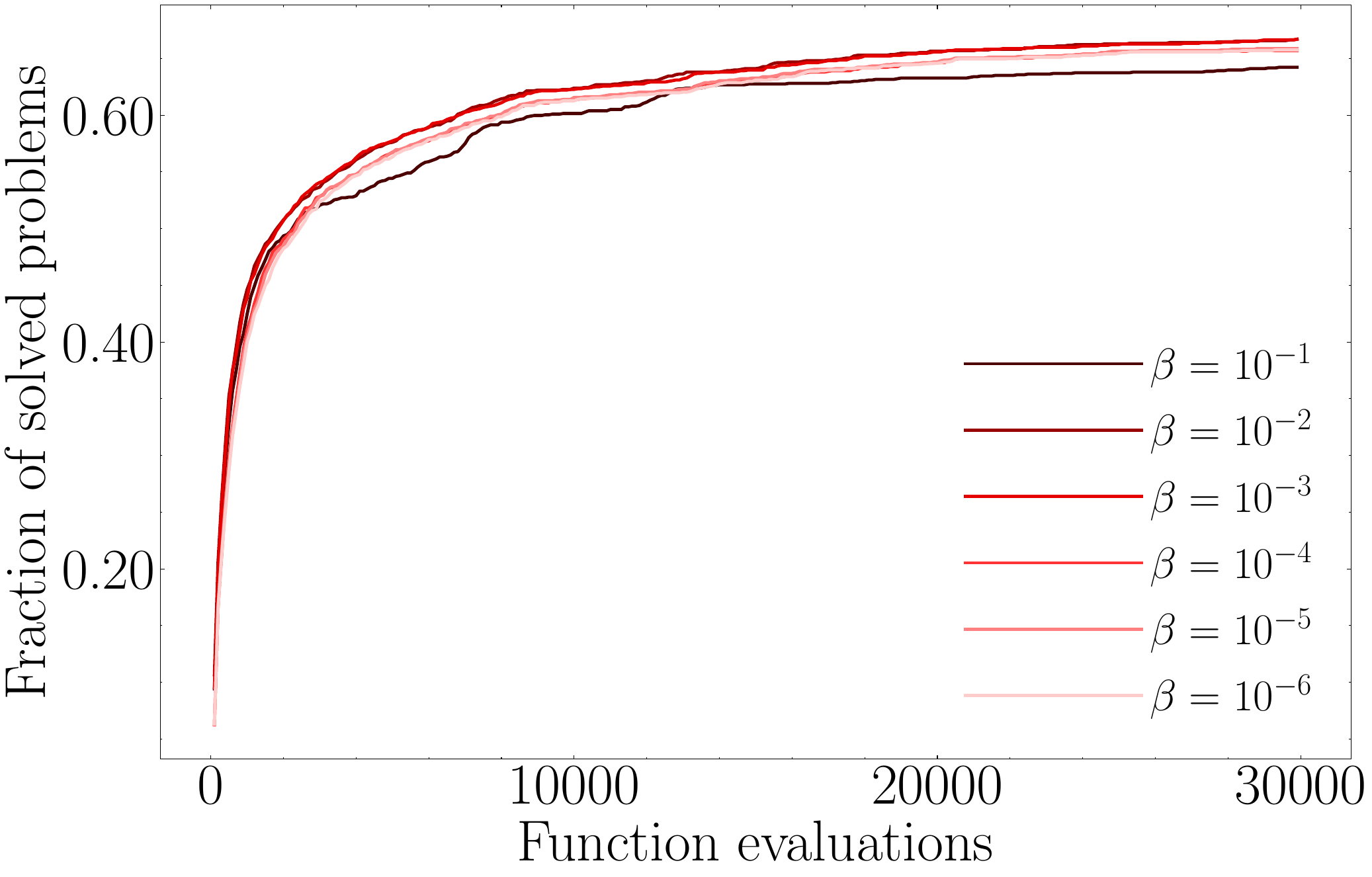}
		\caption{$\text{HALO}_{\text{L-BFGS-B}}$.}
		\label{fig:beta_oc_dfn}
	\end{subfigure}
	\begin{subfigure}[t]{0.48\textwidth}
		\includegraphics[width=1.0\textwidth]{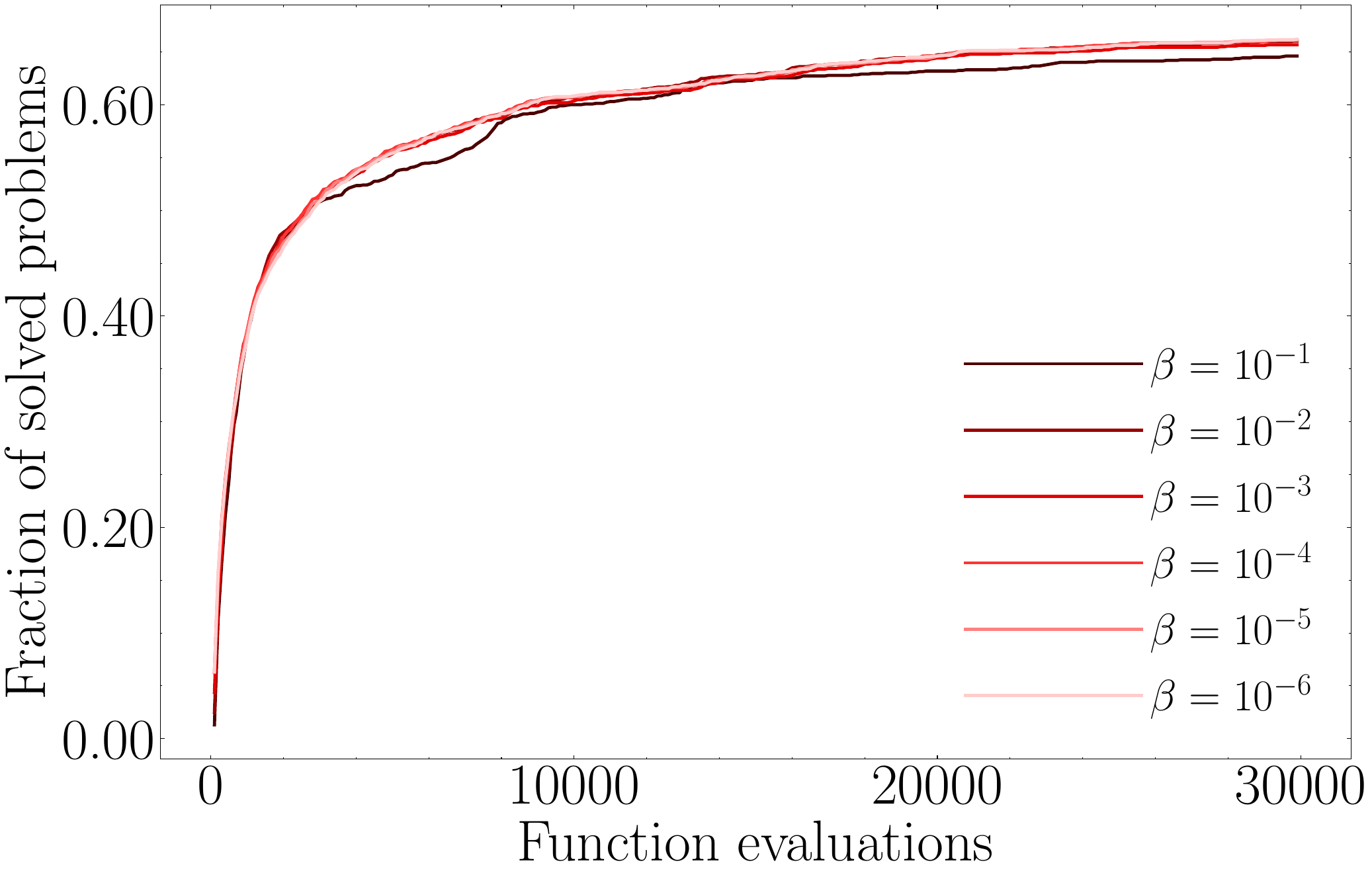}
		\caption{$\text{HALO}_{\text{SDBOX}}$.}
		\label{fig:beta_oc_df}
	\end{subfigure}
	\caption{Operational characteristics for the two versions of HALO varying the parameter $\beta$ regarding the whole set of test functions.}
\end{figure}

In this section, we will analyze the impact of different choices for the parameter $\beta$ on the performance of HALO. 
The parameter $\beta$ is used in Eq. \ref{eq:local_beta} to determine whether a partition $\mathcal{D}_{i_k}$ is small enough to consider its centroid as a starting point for the local optimizer.
A small value of $\beta$ will result in HALO using the local optimizer infrequently, while a higher value will make HALO rely more heavily on the local optimizer.
However, it is important to be cautious when selecting a high value for $\beta$ because once the local search starts from the centroid of a partition, that partition will never be further partitioned in the future as described in Section \ref{section:coupling strategy}.
We conducted experiments using $\text{HALO}_{\text{L-BFGS-B}}$ and $\text{HALO}_{\text{SDBOX}}$ with various values of $\beta$ in the range
\begin{equation}
\beta \in \{10^{-1}, 10^{-2}, 10^{-3}, 10^{-4}, 10^{-5}, 10^{-6}\}
\end{equation}
for all the functions in the three different benchmarks. The results are summarized in Table \ref{table:res_beta}.
\begin{table}[htb!]
	\centering
 \caption{Resume of the numerical results for the sensitivity analysis concerning the parameter $\beta$ for the whole set of test functions.}
	 \resizebox{\textwidth}{!}{\begin{minipage}{\textwidth}
       \begin{tabular}{*7c}\toprule
		\multicolumn{7}{c}{Numerical results} \\
		\midrule
		 & \multicolumn{2}{c}{Percentage solved} &  \multicolumn{2}{c}{Average functions evaluations} & \multicolumn{2}{c}{Average $n^{\circ}$ local searches} \\ 
		$\beta$ &$\text{HALO}_{\text{L-BFGS-B}}$ & $\text{HALO}_{\text{SDBOX}}$ &  
		$\text{HALO}_{\text{L-BFGS-B}}$ & $\text{HALO}_{\text{SDBOX}}$  &  
		$\text{HALO}_{\text{L-BFGS-B}}$ & $\text{HALO}_{\text{SDBOX}}$\\
		\midrule
	    $10^{-1}$ & $64.2$ & $64.6$ & $2230.2$ $(4298.6)$& $2728.3$ $(4836.9)$& $4.5$ $(13.1)$ & $4.6$ $(16.3)$\\
	    $10^{-2}$ & $66.6$ & $66.0$ & $2263.6$ $(4502.2)$& $2744.1$ $(4850.2)$& $2.5$ $(13.5)$ & $2.5$ $(11.4)$\\
	    $10^{-3}$ & $66.6$ & $65.6$ & $2326.9$ $(4660.2)$& $2738.2$ $(4898.8)$ & $1.3$ $(1.0)$ & $1.4$ $(1.0)$\\
		$10^{-4}$ & $65.6$  & $66.0$ & $2403.6$ $(4504.2)$& $2726.8$ $(4891.7)$& $1.1$ $(0.9)$ & $1.1$ $(0.9)$\\ 
		$10^{-5}$ & $65.8$  & $66.1$ & $2426.2$ $(4493.0)$& $2784.0$ $(4988.8)$& $1.0$ $(0.9)$ & $1.0$ $(0.9)$\\ 
		$10^{-6}$ & $65.7$  & $66.2$ & $2472.0$ $(4532.7)$& $2827.2$ $(5039.3)$& $1.0$ $(0.9)$ & $1.0$ $(0.8)$\\
		\bottomrule
	\end{tabular}\end{minipage}}
	\label{table:res_beta}
\end{table}
In general, the performance of HALO appears to be remarkably robust to the choice of the parameter $\beta$, particularly for $\beta$ values between $10^{-2}$ and $10^{-6}$. This observation is supported by the operational characteristics shown in Fig. \ref{fig:beta_oc_dfn} and Fig. \ref{fig:beta_oc_df} for $\text{HALO}_{\text{L-BFGS-B}}$ and $\text{HALO}_{\text{SDBOX}}$, respectively.
Moreover, within the range of $\beta$ values from $10^{-2}$ to $10^{-6}$, the choice of the local optimizer (L-BFGS-B or SDBOX) does not significantly affect the performance of HALO. Although $\text{HALO}_{\text{L-BFGS-B}}$ demonstrates the best performance in terms of average function evaluations and percentage of problems solved, the difference in performance between $\text{HALO}_{\text{SDBOX}}$ and $\text{HALO}_{\text{L-BFGS-B}}$ is not substantial.

\begin{figure}[htb!]
\centering
	\begin{subfigure}[t]{0.29\textwidth}
		\includegraphics[width=1\textwidth]{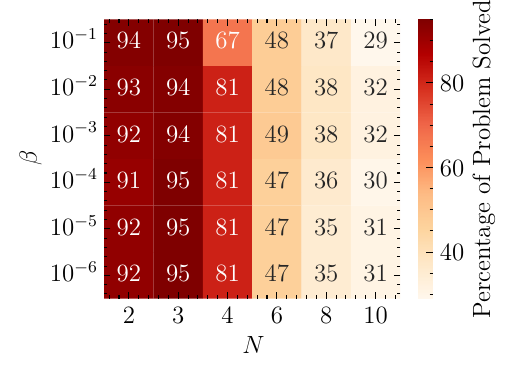}
		\caption{$\text{HALO}_{\text{L-BFGS-B}}$.}
		\label{fig:dim_beta_dfn}
	\end{subfigure}
	\begin{subfigure}[t]{0.29\textwidth}
    	\includegraphics[width=1\textwidth]{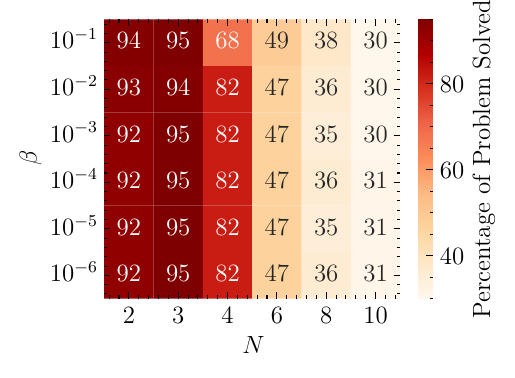}
    	\caption{$\text{HALO}_{\text{SDBOX}}$.}
    	\label{fig:dim_beta_df}
	\end{subfigure}
	\caption{Percentage of problems solved varying the dimensionality $N$ conditional to the value of $\beta$ for the two versions of HALO.}
	\label{fig:dim_beta}
\end{figure}

Table \ref{table:res_beta} also presents the average number of local optimization routines initiated varying $\beta$. On average, only one or two local optimization routines start when $\beta \in [10^{-3}, 10^{-6}]$ or when $\beta = 10^{-2}$, respectively, for both $\text{HALO}_{\text{SDBOX}}$ and $\text{HALO}_{\text{L-BFGS-B}}$. 
However, the situation changes when $\beta = 10^{-1}$, as both $\text{HALO}_{\text{SDBOX}}$ and $\text{HALO}_{\text{L-BFGS-B}}$ experience a performance decline, especially for $N=4$ dimensions (see Fig. \ref{fig:dim_beta}).

Based on the results, we recommend the following regarding the choice of $\beta$:
\begin{itemize}
    \item Avoid selecting a value of $\beta$ less than $10^{-4}$ to prevent HALO from consuming excessive function evaluations on very small hyperrectangles and to avoid potential numerical instabilities.
    \item If the objective function is expected to contain noise, we recommend using the derivative-free version of HALO, namely $\text{HALO}_{\text{SDBOX}}$, with $\beta$ values between $10^{-2}$ and $10^{-4}$.
    Otherwise, $\text{HALO}_{\text{L-BFGS-B}}$ should be the default choice, again with $\beta$ values between $10^{-2}$ and $10^{-4}$.
\end{itemize}
\subsection{Summary of the Numerical Results and Future Work}
Based on the numerical results from the benchmark, we will briefly summarize the pros and cons of HALO. The numerical experiments highlight the following advantages of HALO:
\begin{itemize}
   \item[\textcolor{green}{\checkmark}] \textbf{Fast identification of the global minimum}: HALO shows a strong ability to rapidly converge to the global minimum, balancing exploitation and exploration and significantly outperforming many popular DIRECT-type and stochastic methods in terms of convergence speed. Numerical results highlight that our proposed estimates of the local Lipschitz constants are crucial for achieving this high performance.
   \item[\textcolor{green}{\checkmark}] \textbf{Highlights problem interpretability}: HALO provides valuable insights into the optimization process by identifying the importance of different variables (see Section \ref{subsec:varimp}), which is crucial for understanding complex objective functions and informing decision-making. This interpretability can also facilitate dimensionality reduction by allowing less significant variables to be excluded from the analysis, simplifying the problem-solving process.
   \item[\textcolor{green}{\checkmark}] \textbf{Hyperparameter free}: HALO operates without the need for critical hyperparameters that significantly affect convergence, unlike similar methods that rely on estimating local Lipschitz constants. Additionally, our hyperparameter sensitivity analysis in Section \ref{subsec:PSA} demonstrates the remarkable robustness of HALO regarding the hyperparameter that controls the number of local search initializations.
   \item[\textcolor{green}{\checkmark}] \textbf{Insensitive to the choice of the local optimizer}: 
    We tested HALO with two local optimizers: a quasi-Newton method and a derivative-free method. Numerical results showed no significant difference in performance when varying the local optimizers. We attribute this to the local search being initiated only in the proximity of a stationary point for solution refinement, limiting the potential for varying behavior. This flexibility allows HALO to handle many types of global optimization problems, whether they are differentiable or not.  
\end{itemize}
As with any algorithm also HALO has drawbacks that can be summarized:
\begin{itemize}
   \item[\textcolor{red}{$\times$}] \textbf{Computational complexity}: We did not compare algorithms based on computational time because, in the context of black-box optimization, the objective function typically originates from a simulation process that constitutes the majority of the computational burden. Function evaluations are often limited, making the computational complexity of the optimization algorithm itself negligible. However, it's important to note that HALO performs additional operations, such as evaluating local Lipschitz constants, which increases its computational complexity compared to some simpler algorithms. This added complexity can become a limitation when the number of allowed function evaluations is very high (greater than 50,000).
   \item[\textcolor{red}{$\times$}] \textbf{Sensitive to the problem dimension}: Although HALO has demonstrated strong performance across various benchmarks, like any global optimization algorithm, its effectiveness can diminish as the problem dimensionality increases due to the 'curse of the dimensionality' \cite{bellman1966dynamic}. The current version of HALO is therefore best suited for problems of low to medium dimensionality.
\end{itemize}
Given those pros and cons we can identify the following potential directions and guidelines for future work and development:
\begin{itemize}
     \item\textbf{Dimensionality Reduction}: A promising direction for future research is to use the variable importance derived from HALO to guide dimensionality reduction, which is especially relevant for real-world applications \cite{d2024generative}. One approach would be to leverage the approximated gradients to define an active subspace \cite{constantine2014active, d2024learning}, thus focusing on the most influential dimensions of the problem. This can potentially reduce computational costs, make the optimization more efficient, and speed up convergence to the global minimum.
   \item\textbf{Improving the hybridization strategy}: HALO currently employs the local optimization algorithm infrequently within its hybrid scheme. Increasing the utilization of local optimization could enhance HALO's performance, particularly in high-dimensional spaces. This improvement might also reduce the algorithm's computational time.
    \item\textbf{Considering other partition schemes}: Our algorithm uses a central partition scheme to divide the search space but other techniques, such as those proposed in  \cite{paulavivcius2014simplicial}, \cite{stripinis2024lipschitz} could also be considered.
   \item\textbf{Better convergence properties}: While HALO leverages the everywhere dense property to guarantee a path to the global minimum, future refinements may further strengthen its convergence guarantees or improve stopping criteria. Such improvements should preserve the algorithm’s practical usability and avoid introducing new, critical hyperparameters.. 
\end{itemize}
\section{Conclusions}
In this study, we introduced HALO (Hybrid Adaptive Lipschitzian Optimization), a deterministic partition-based global optimization algorithm. HALO incorporates an adaptive procedure to estimate the local Lipschitz constant for each partition or sub-region $\mathcal{D}_i$ within the box domain $\mathcal{D}$. The estimation is performed automatically based on the size of the partition, eliminating the need for the user to define crucial hyperparameters in advance. The adaptive procedure strikes a balance between the current estimate of the global Lipschitz constant and the norm of the gradient approximation associated with the centroids of each partition.

We proposed a simple coupling strategy with local optimization algorithms, including one gradient-based and one derivative-free, to expedite convergence towards a stationary point. By leveraging these local optimization techniques, HALO achieves faster convergence speeds, enhancing its effectiveness in solving optimization problems. This flexibility allows HALO to function effectively as a derivative-free algorithm when paired with a method like SD-BOX. However, recognizing that some applications may benefit from gradient information, HALO is also adaptable. Users can opt for a quasi-Newton scheme for local search if accurate derivative computation is feasible. This dual capability ensures HALO can handle a wide range of global optimization problems, both differentiable and non-differentiable.

To evaluate the performance of HALO, we conducted a comprehensive comparison with popular GO algorithms using three diverse benchmarks comprising numerous test functions. The results demonstrated that HALO exhibited robust performance across different problem domains. Additionally, we conducted a sensitivity analysis on a hyperparameter controlling the local search, and the findings indicated that HALO's performance remains consistent within a certain range of hyperparameter values. 

HALO not only enables efficient global optimization but also facilitates the extraction of valuable insights from the objective function. By analyzing the matrix that collects the approximation of the gradients around each partition's centroid, we can identify the influential directions in the search space and gain a deeper understanding of the factors driving the reduction of the objective function. This feature makes HALO a powerful tool not only for optimization but also for exploratory analysis and problem understanding in a wide range of applications.

Based on these results, HALO emerges as highly competitive and can extend our tools of methodologies for tackling complex real-world GO problems. Its deterministic nature, lack of crucial hyperparameters to set, and integration with gradient or derivative-free local optimization methods position HALO as a valuable addition to the toolbox of optimization practitioners.
\FloatBarrier

\begin{figure}[htb!]
\centering
	\begin{subfigure}[t]{0.29\textwidth}
		\includegraphics[width=1\textwidth]{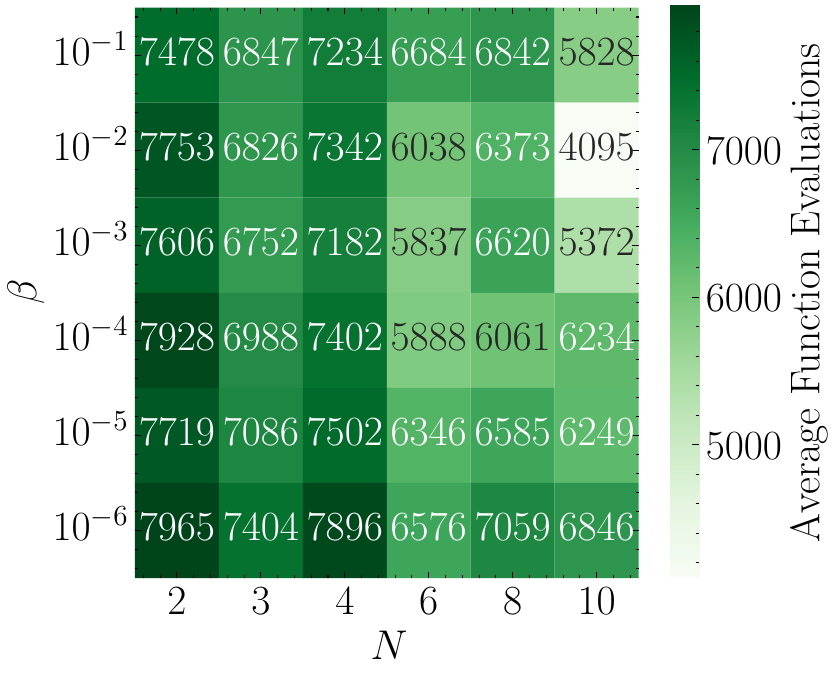}
		\caption{$\text{HALO}_{\text{L-BFGS-B}}$.}
		\label{fig:dim_beta_fe_dfn}
	\end{subfigure}
	\begin{subfigure}[t]{0.29\textwidth}
    	\includegraphics[width=1\textwidth]{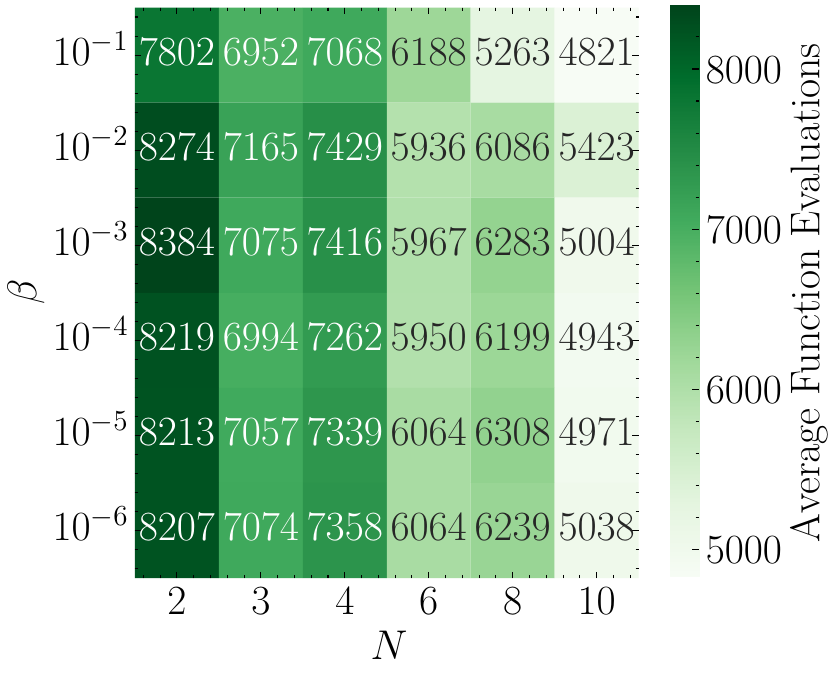}
    	\caption{$\text{HALO}_{\text{SDBOX}}$.}
    	\label{fig:dim_beta_fe_df}
	\end{subfigure}
	\caption{Average number of functions evaluations varying the dimensionality $N$ conditional to the value of $\beta$ for the two versions of HALO.}
	\label{fig:dim_beta_fe}
\end{figure}
\subsection*{Data Availability Statement}
Data sharing not applicable to this article as no datasets were generated or analysed during the current study.
\subsection*{Acknowledgment}
The author would like to thank Prof. Stefano Lucidi for his valuable advice and guidance throughout this work. Appreciation is also extended to the three anonymous reviewers for their constructive comments, which helped improve the quality of the paper.
\bibliography{biblio}
\end{document}